\documentclass[12pt]{amsart}
\usepackage{amsmath}  
\usepackage{amsthm}
\usepackage{amsfonts}
\usepackage{amssymb}
\usepackage{amstext}
\usepackage{amscd}
\usepackage{enumerate}
\usepackage[all]{xy}
\newtheorem{thm}{Theorem}[section]
\newtheorem{prop}[thm]{Proposition}
\newtheorem{cor}[thm]{Corollary}
\newtheorem{lem}[thm]{Lemma}

\theoremstyle{definition}

\newtheorem{dfn}[thm]{Definition}
\newtheorem{rem}[thm]{Remark}

\newtheorem{conj}[thm]{Conjecture}

\def\Z{\mathbb{Z}}
\def\Q{\mathbb{Q}}
\def\R{\mathbb{R}}
\def\C{\mathbb{C}}

\def\E{\mathbb{E}}
\def\cA{\mathcal{A}}

\def\cG{\mathcal{G}}
\def\cL{\mathcal{L}}
\def\fa{\mathfrak{a}}
\def\fb{\mathfrak{b}}
\def\fc{\mathfrak{c}}
\def\fg{\mathfrak{g}}

\def\fs{\mathfrak{s}}
\def\su{\mathfrak{su}}
\def\ad{\operatorname{ad}}
\def\conj{\operatorname{conj}}

\def\im{\operatorname{im}}
\def\Image{\operatorname{Im}}
\def\ind{\operatorname{ind}}
\def\Ind{\operatorname{Ind}}
\def\inst{\operatorname{inst}}
\def\Ker{\operatorname{Ker}}
\def\Coker{\operatorname{Coker}}
\def\Hom{\operatorname{Hom}}
\def\sign{\operatorname{sign}}
\def\Tr{\operatorname{Tr}}
\def\tM{\tilde{M}}
\def\M{M_{\rho \sigma}}
\def\Dir{\not \! \! \mathfrak{D}}
\def\bd{\partial^{(0)}}
\def\cm{C^{(0)}}
\def\bdg{\partial_{\gamma}}
\def\rs{\rho \sigma}

\def\CP{\mathbb{CP}}
\def\dbar{\bar{\partial}}
\begin{document}

\title[Instanton Floer homology for lens spaces]
{Instanton Floer homology for lens spaces}
\author{Hirofumi Sasahira}
\date{}
\renewcommand{\thefootnote}{\fnsymbol{footnote}}
\footnote[0]{2000\textit{ Mathematics Subject Classification}.
Primary 57R57, 57R58.}
\address
{Graduate school of Mathematics, Nagoya University, \endgraf
Furocho,  Chikusaku, Nagoya, Japan.  464-8602          
}
\email{hsasahira@math.nagoya-u.ac.jp}

\maketitle

\begin{abstract}
Floer constructed instanton homology for homology 3-spheres. In this paper, we extend instanton Floer homology to lens spaces $L(p,q)$. Moreover we show a gluing formula for a variant of Donaldson invariant along lens spaces. As an application, we prove that $X = \CP^2 \# \CP^2$ does not admit a decomposition $X = X_1 \cup X_2$. Here $X_1$ and $X_2$ are oriented, simply connected, non-spin 4-manifolds with $b^+ = 1$ and with boundary $L(p, 2)$, and $p$ is a prime number of the form $16N + 1$.
\end{abstract}

\section{Introduction}
\label{intro}
Instanton Floer homology $HF_*(Y)$ was constructed in \cite{Fl} for oriented homology 3-spheres $Y$. This invariant is defined using flat connections over $Y$ and moduli spaces of instantons over $Y \times \R$.
 As is well known, instanton Floer homology has an important role in calculations of Donaldson invariants for closed 4-manifolds. Let $X$ be an oriented closed 4-manifold with $b^+(X) > 1$ and fix a cohomology class $c \in H^2(X;\Z)$. Donaldson invariant $\Psi_{X,c}$ for $X$ is defined as a $\Q$-valued function on $A(X) = \oplus_{ d \geq 0 } H_2(X;\Z)^{ \otimes d}$, using moduli spaces $M_P$ of instantons on principal $U(2)$-bundles $P$ with $c_1(P) = c$.  Suppose that $X$ has a decomposition $X = X_1 \cup X_2$, where $X_1$ and $X_2$ are compact 4-manifolds with $b^+ > 0$ and with boundary $Y$ and $-Y$ respectively. Here $-Y$ is $Y$ with opposite orientation. We can define relative Donaldson invariants $\Psi_{X_1,c_1}:A(X_1) = \otimes H_2(X_1;\Z) \rightarrow HF_*(Y)$, $\Psi_{X_2,c_2}:A(X_2) =\otimes H_2(X_2;\Z) \rightarrow HF_*(-Y)$ where $c_1 = c|_{X_1}, c_2 = c|_{X_2}$.
There is a natural pairing $< \cdot , \cdot >:HF_*(Y) \otimes HF_*(-Y) \rightarrow \Q$, and  we have a gluing formula $\Psi_{X,c} = < \Psi_{X_1, c_1}, \Psi_{X_2,c_2} >$. Note that $A(X) = A(X_1) \otimes A(X_2)$ since $Y$ is a homology 3-sphere. We can completely  determine $\Psi_{X,c}$ in terms of the relative invariants from the gluing formula. 

There is a variant $\Psi_{X}^{u_1}$ of Donaldson invariants \cite{FS}, \cite{S torsion} defined using a cohomology class $u_1 \in H^1(M_P;\Z_2)$. The cohomology class $u_1$ is the first Stiefel-Whitney class of the determinant line bundle of the real part of twisted Dirac operators over $X$.
The variant $\Psi_{X}^{u_1}$ is a function on a subspace $A'(X)$ of $A(X)$ with values in $\Z_2$. Variants of instanton Floer homology are defined for oriented homology 3-spheres in \cite{FFO} and \cite{S Floer}, and there is a similar gluing formula for $\Psi_X^{u_1}$.

A natural problem is how to define instanton Floer homology for more general 3-manifolds $Y$, which enable us to construct gluing formulas for more general decompositions of $X$.  Mainly there are two difficulties when we try to generalize instanton Floer homology:

\begin{enumerate}[(i)]
\item
The existence of {\it reducible} (projectively) flat connections on 3-manifolds.

\vspace{2mm}

\item
$H_2(X;\Z)$ is {\it not} isomorphic to the direct sum $H_2(X_1;\Z) \oplus H_2(X_2;\Z)$.

\end{enumerate}

(i) implies that moduli spaces of instantons over $Y \times \R$ can be singular. Even if the moduli spaces are smooth, the usual proof that the square of the boundary map is zero breaks down. (ii) implies that we must consider the situation where a surface in $X$ representing an element in $H_2(X;\Z)$ is also decomposed by $Y$.

 There are some partial answers to each problem. Austin-Braam \cite{AB} and Donaldson \cite{D Floer} introduced equivalent versions of instanton Floer homology under some assumptions in order to overcome the problem (i) of reducible flat connections. The equivalent Floer homologies enable us to generalize the gluing formula for Donaldson invariants $\Psi_{X,c}$. Furuta \cite{Fur inv inst} defined an analog $I_*(L(p,q))$ of instanton Floer homology for lens spaces $L(p,q)$, making use of Dirac operators over $L(p,q) \times \R$. Note that all flat connections on $L(p,q)$ are reducible since the fundamental group of $L(p,q)$ is abelian. Using this analog we can construct a gluing formula for the variant $\Psi_X^{u_1}$ of Donaldson invariants.  However these gluing formulas can not be applied to the problem (ii). That is, these gluing formulas calculate only the restrictions of the invariants $\Psi_{X,c}$, $\Psi_{X}^{u_1}$ to the images of the natural maps $A(X_1) \otimes A(X_2) \rightarrow A(X)$, $A'(X_1) \otimes A'(X_2) \rightarrow A'(X)$. 
 
On the other hand, Fukaya \cite{Fuk} introduced a generalization $HFF_*(Y, Q; \gamma)$ of instanton Floer homology for general 3-manifolds $Y$, $U(2)$-bundles $Q$ over $Y$ and a loop $\gamma$ in $Y$, provided that all (projectively) flat connections on $Q$ are irreducible. This generalization gives a complete answer to (ii). That is, we can show a gluing formula which completely calculate Donaldson invariants in terms of relative invariants of $X_1, X_2$. See \cite{BD}. But it seems that the assumption on (projectively) flat connections on $Q$ has not been removed. Thus the problems (i) and (ii) have been separately dealt with.
 
In this paper, we deal with both (i) and (ii) at the same time for lens spaces. We apply Fukaya-Floer type construction to lens spaces with some modifications (Subsection \ref{ss FF}). We will define an analog $I_*(L(p,q);\gamma)$ of Fukaya-Floer homology for an odd prime integer $p$ and a loop $\gamma$ in $L(p,q)$. In the construction, we make use of Dirac operators as in \cite{Fur inv inst}. Moreover we construct a gluing formula for $\Psi_{X}^{u_1}$ along $L(p,q)$.

As an application, we will prove that $X = \CP^2 \# \CP^2$ does not admit a decomposition $X = X_1 \cup X_2$. Here $X_1$ and $X_2$ are simply connected, non-spin 4-manifolds with $b^+ = 1$ and with boundary $L(p,2)$ and $-L(p,2)$ respectively, and $p$ is a prime number of the form $16N+1$. See Theorem \ref{thm dec}.
This is based on a calculation  of $I_*(L(p,2)), I_*(L(p,2);\gamma)$ and the non-vanishing of $\Psi_X^{u_1}$. The calculation of Floer homologies requires counting the number of instantons over $L(p,2) \times \R$. This was done in \cite{Au}, \cite{FH} and \cite{Fur inv inst}. The non-vanishing of $\Psi_X^{u_1}$ was proved in \cite{S Floer}.

We give a remark which is related to Seiberg-Witten theory. In \cite{W}, Witten introduced Seiberg-Witten equations and defined Seiberg-Witten invariants using the moduli spaces of solutions to the equations. Witten also conjectured that Seiberg-Witten invariants are equivalent to Donaldson invariants and that Donaldson invariants can be calculated in terms of Seiberg-Witten invariants through a formula. This formula has been proved for many 4-manifolds. (See \cite{FL} and \cite{GNY}.) Moreover Seiberg-Witten theory gives us simpler proofs of many results obtained by Donaldson theory and new stronger results.
Theorem \ref{thm dec} is in contrast to such things.
We should not expect that Theorem \ref{thm dec} can be proved by Seiberg-Witten theory, because $\CP^2 \# \CP^2$ has a metric of positive scalar curvature and any invariants from Seiberg-Witten equations (Seiberg-Witten invariants and a refinement due to Bauer and Furuta \cite{BF}) are trivial.

\section{Constructions of instanton homology}

\subsection{Preliminaries}
Let $p, q$ be relatively prime integers with $0 < q < p$, and denote by $L(p, q)$ the lens space $S^3/\Z_p$. Here the action of $\Z_p = \{ \ \zeta \in \C \ | \ \zeta^p = 1 \ \}$ on $S^3 = \{ \ (z_1, z_2) \in \C^2 \ | \ |z_1|^2 + |z_2|^2 = 1 \ \}$ is defined by
\[
\zeta \cdot ( z_1, z_2) = (\zeta z_1, \zeta^q z_2).
\]
Throughout this paper, we consider only the Riemannian metric on $L(p,q)$ induced by the standard Riemannian metric on $S^3$.
In this subsection, we take up some basic facts about $SU(2)$-flat connections over $L(p,q)$ and moduli spaces of instantons over $L(p,q) \times \R$.
\vspace{2mm}

We consider $L(p,q)$ as an oriented manifold with the orientation induced by the standard orientation of $S^3$. We write $-L(p,q)$ for $L(p,q)$ with the opposite orientation. Let $Y$ be $L(p,q)$ or $-L(p,q)$. The moduli space $R(Y)$ of flat connections on the trivial $SU(2)$-bundle $Q = Y \times SU(2)$ is identified with $\Hom (\pi_1(Y), SU(2))/\conj$. Since $\pi_1(Y)$ is abelian, all flat connections are reducible. That is, the stabilizer $\Gamma_{\rho}$ of any flat connection $\rho$ in the gauge group is isomorphic to $U(1)$ or $SU(2)$.
For each class $[\rho] \in R(Y)$ represented by a flat connection $\rho$, we define an index $\delta_Y([\rho]) \in \Z_8 = \Z/ 8 \Z$ as follows. Let $A$ be an $SU(2)$-connection over $Y \times \R$ such that 
\[
A = \left\{
\begin{array}{cl}
\pi^* \rho & \text{on $Y \times (-\infty, -1)$,} \\
\pi^* \theta & \text{on $Y \times (1, \infty)$.} 
\end{array}
\right.
\]
Here $\pi:Y \times \R \rightarrow Y$ is the projection and $\theta$ is the trivial flat connection of $Q$. Take a small positive number $\epsilon > 0$ and define the function $W^{-+}:Y \times \R \rightarrow \R_{>0}$ by
\[
W^{-+}(y, t) = e^{\epsilon t}.
\]
We define a weighted $L^2$ norm $\| \cdot \|_{L^{2,(-\epsilon, \epsilon)}}$ on the sections of $(\Lambda^0_{Y \times \R} \oplus \Lambda^+_{Y \times \R}) \otimes \pi^* \fg_{P}$ by
\[
\| f \|_{L^{2, (-\epsilon, \epsilon)}} := \| W^{-+} f \|_{L^2}.
\]

Similarly, we define a weighted $L^2_1$ norm $\| \cdot \|_{L_1^{2, (-\epsilon, \epsilon)}}$ on the sections of $\Lambda^1_{Y \times \R} \otimes \pi^* \fg_P$ by
\[
\| f \|_{L^{2, (-\epsilon, \epsilon)}_1} := \| W^{- +} f \|_{L^2_1}.
\]
We consider the operator
\[
D_A=d_A^* + d_A^+:
L^{2,(-\epsilon, \epsilon)}_1 (\Lambda_{Y \times \R}^1 \otimes \pi^*\fg_{P})
\longrightarrow
L^{2, (-\epsilon, \epsilon)}( (\Lambda^0_{Y \times \R} \oplus \Lambda^+_{Y \times \R}) \otimes \fg_P ).
\]
We define $\delta_Y([\rho]) \in \Z_8$ to be $\ind^{-+} D_A \mod 8$. Here $\ind^{-+} D_A := \dim \Ker D_A - \dim \Coker D_A$. We can easily see that $\delta_{Y}([\rho])$ depends only on the class $[\rho]$.

\vspace{3mm}

For each flat connection $\rho$ over $Y$, we have the complex:
\[
\Omega^0_Y(\fg_P) \stackrel{d_{\rho}}{\longrightarrow}
\Omega^1_Y(\fg_P) \stackrel{d_{\rho}}{\longrightarrow}
\Omega^2_Y(\fg_P).
\]
Let $H^i(Y;\ad \rho)$ be the $i$-th cohomology group of this complex.

\begin{lem} \label{lem H^1}
$H^1(Y;\ad \rho)$ is trivial.
\end{lem}

To prove this, consider the pull-back $\tilde{\rho}$ of $\rho$ by the projection $S^3 \rightarrow Y=S^3/\Z_p$. Then $H^1(Y;\ad \rho)$ is identified with the invariant subspace of the natural action of $\Z_p$ on $H^1(S^3; \ad \tilde{\rho})$. Since $S^3$ is simply connected, $\tilde{\rho}$ is gauge equivalent to the trivial connection. This means that $H^1(S^3;\ad \tilde{\rho})$ is isomorphic to $H^1(S^3;\R) \otimes \su(2)$. But $H^1(S^3;\R)$ is trivial. Hence we have obtained the statement.

\vspace{2mm}

Lemma \ref{lem H^1} implies that if the curvature of an instanton over $Y \times \R$ is $L^2$-integrable, then the instanton exponentially converges to some flat connections at $\pm \infty$ with respect to any Sobolev norms. (See \cite[Section 4.1]{D Floer}.) We consider moduli spaces of instantons whose curvatures are $L^2$-integrable. Let $\tilde{M}_{\rho \sigma}$ denote the framed moduli space of instantons with limits $\rho, \sigma$. That is, $\tilde{M}_{\rho \sigma}$ is the quotient of the space of instantons with limits $\rho$ at $-\infty$ and $\sigma$ at $+\infty$ by the group of gauge transformations with limit $1$ at $\pm \infty$. The group $\Gamma_{\rho} \times \Gamma_{\sigma}$ naturally acts on $\tilde{M}_{\rho \sigma}$, and put $M_{\rho \sigma} := \tilde{M}_{\rho \sigma} / \Gamma_{\rho} \times \Gamma_{\sigma}$. ( As stated above, $\Gamma_{\rho}$ is the stabilizer of $\rho$ in the gauge group.)

\begin{lem} \label{lem 1}
Let $\rho, \sigma$ be flat connections over $Y$ which represent deferent classes in $R(Y)$. Then the moduli space $M_{\rho \sigma}$ is a smooth manifold, and
\[
\dim M_{\rho \sigma} \equiv 
\delta_{Y}([\rho]) - \delta_Y([\sigma]) - \dim \Gamma_{\rho} \mod 8.
\]
\end{lem}

First we consider the deformation complex of the framed moduli space $\tilde{M}_{\rho \sigma}$ at $[A]$:
\begin{equation} \label{eq def cpx fr}
L^{2, (\epsilon,\epsilon)}_2(\Lambda^0_{Y \times \R} \otimes \su(2))
\stackrel{d_A}{\longrightarrow}
L^{2, (\epsilon,\epsilon)}_1(\Lambda^1_{Y \times \R} \otimes \su(2))
\stackrel{d_A^+}{\longrightarrow}
L^{2, (\epsilon,\epsilon)}(\Lambda^+_{Y \times \R} \otimes \su(2)).
\end{equation}
Here $L^{2, (\epsilon, \epsilon)}$ is the completion of the space of compact supported sections by a weighted $L^2$ norm defined by a function $W^{++}:Y \times \R \rightarrow \R_{>0}$ with
\[
W^{++}(y,t) = \left\{
\begin{array}{ll}
e^{- \epsilon t} & \text{ if $t < -1$} \\
e^{\epsilon t} & \text{ if $t > 1$}.
\end{array}
\right.
\]
Similarly for $L^{2, (\epsilon, \epsilon)}_1, L^{2, (\epsilon, \epsilon)}_2$. We can show that the second cohomology of the complex is trivial by using Weitzenb\"ock formula for $d_A^{*, \epsilon} + d^+_A$, since the Riemannian metric of $Y \times \R$ is self-dual and the scalar curvature is positive. (See \cite{AHS} for the case when the 4-manifold is closed.) Here $d^{*,\epsilon}$ is the formal adjoint of $d_A$ with respect to the wighted Sobolev norms. This implies that $\tilde{M}_{\rho \sigma}$ is a smooth manifold. 

Next we show that all instantons $A$ with limits $\rho, \sigma$ are irreducible. Suppose that $A$ is reducible. Then we can write
\[
A = a \oplus -a
\]
for some $U(1)$-connection $a$. Since $*F_a = -F_a$ and $F_a$ is closed, we have
\[
d^* F_a = * d * F_a = - *dF_a = 0.
\]
Hence $F_a$ is a harmonic 2-form over $Y \times \R$. Moreover $F_a$ decays exponentially as $t \rightarrow \pm \infty$. Thus we have
\[
F_a \in \ker (d^{*} + d) \cap L^2 \cong 
\Image (H^2(X_T, \partial X_T;\R) \rightarrow H^2(X_T;\R)).
\]
Here $T>0$ is a positive number and $X_T = Y \times [-T, T]$. Note that $H^2(X_T;\R) = 0$ since $Y$ is a lens space. This means that $a$ (and hence $A$) is a flat connection over $Y \times \R$. 

Let $A'$ be a connection which is gauge equivalent to $A$ and is in temporal gauge. Then by the instanton equation we have
\[
\frac{ \partial A' }{\partial t} = -*_Y F_{A'} = 0,
\]
where $*_Y$ is Hodge $*$-operator over $Y$ and we have used $F_{A'} = 0$. Therefore the restriction $A'_t$ of $A'$ to $Y \times \{ t \}$ is independent of $t$, and especially $[\rho] = [\sigma]$ in $R(Y)$. This is a contradiction since we assumed $[\rho] \not= [\sigma]$. Thus $A$ is irreducible.

\vspace{2mm}

The fact that $A$ is irreducible implies that the stabilizer of $[A] \in \tilde{M}_{\rho \sigma}$ in $\Gamma_{\rho} \times \Gamma_{\sigma}$ is $\{ \pm (1, 1) \}$ and  that the action of $\Gamma_{\rho} \times \Gamma_{\sigma} / \{ \pm (1,1) \}$ on $\tilde{M}_{\rho \sigma}$ is free. Hence the quotient $M_{\rho \sigma}=\tilde{M}_{\rho \sigma}/\Gamma_{\rho} \times \Gamma_{\sigma}$ is also smooth.

\vspace{2mm}

We show the second part of the lemma. The dimension of the framed moduli space $\tilde{M}_{\rho \sigma}$ is the index $\ind^{++}(d_A^{*,\epsilon} + d_A^+)$ of the complex (\ref{eq def cpx fr}). This is equal to
\[
\ind^{-+}(d_A^* + d_A^+) + \dim \Gamma_{\sigma}.
\]
(See \cite[Proposition 3.10, Proposition 3.19]{D Floer}.) Hence we have
\[
\begin{split}
\dim M_{\rho \sigma} 
&=\dim \tilde{M}_{\rho \sigma} - \dim \Gamma_{\rho} - \dim \Gamma_{\sigma} \\
&=\ind^{-+}(d_A^* + d_A^+) - \dim \Gamma_{\rho} \\
&\equiv \delta_{Y}([\rho])-\delta_{Y}([\sigma]) - \dim \Gamma_{\rho} \mod 8.
\end{split}
\]
Here we used the additivity of the index in the last equality.

\subsection{Analog of Floer homology}

In this subsection, we review the construction in \cite{Fur inv inst}. Floer homology for a homology 3-sphere $Z$ \cite{Fl} is defined to be the homology of the chain complex generated by gauge equivalence classes of flat connections over $Z$. The boundary operator is defined by counting number of points of 0-dimensional moduli spaces of instantons over $Z \times \R$ (with signs). If we apply this construction to a lens space $Y$, as explained in \cite{Fur inv inst}, we will face the problem that the square of the boundary operator is {\it not} zero. The idea to overcome this problem is that we modify the definition of the boundary operator using twisted Dirac operators over $Y \times \R$. 

\vspace{3mm}

For each $i \in \Z$, let $CF_i(Y)$ be the vector space over $\Z_2$ spanned by
\[
\{ \ [\rho] \in R(Y) \ | 
\ \Gamma_{\rho} \cong U(1), \ \delta_Y([\rho]) \equiv i \mod 8 \ \}.
\]
Then we put
\[
C^{(0)}_i(Y) := CF_{2i}(Y), \quad
C^{(1)}_i(Y) := CF_{2i+1}(Y).
\]
We will define the boundary operator $\partial^{(0)}:C_*^{(0)} \rightarrow C_{*-1}^{(0)}$ as follows. ( The definition of $\partial^{(1)}:C^{(1)}_* \rightarrow C^{(1)}_{*-1}$ is similar.)

Take generators $[\rho] \in C_i^{(0)}(Y), [\sigma] \in C_{i-1}^{(0)}(Y)$. By the dimension formula in Lemma \ref{lem 1}, we have
\[
\dim M_{\rho \sigma} \equiv 2 - 1 \equiv 1 \mod 8.
\]
We can take representatives $\rho, \sigma$ of the classes such that
\[
\dim M_{\rho \sigma} = 1.
\]
We define $M_{\rho \sigma}' \subset M_{\rho \sigma}$ to be the moduli space of instantons with center of mass $0$. Here the center of mass of $A$ is defined by
\[
\int_{Y \times \R} t | F_A |^2 d \mu_{Y \times \R} \in \R.
\]
Standard arguments, which can be found in \cite{D Floer}, show that $M_{\rho \sigma}'$ is a compact smooth manifold of dimension $0$. That is, $M_{\rho \sigma}'$ is a finite set. 

Fix a spin structure $\fs$ of $Y$ and a connection $A$ with limit $\rho, \sigma$. Then we have the twisted Dirac operator over $Y \times \R$:
\[
\Dir_{A}:
L^{2,(-\epsilon, \epsilon)}_1( S^+ \otimes E) \longrightarrow 
L^{2,(-\epsilon, \epsilon)} (S^- \otimes E).
\]
Here $E$ is the rank $2$ complex vector bundle over $Y \times \R$ associated with $\pi^* Q$ and $S^{\pm}$ are the spinor bundle of the spin structure. We denote $\ind^{-+} \Dir_A \in \Z$ by $i_{\rho \sigma}$.
We put
\[
< \partial^{(0)}([\rho]), [\sigma] > := 
\left\{
\begin{array}{cll}
\# \M'  &\mod 2 & \text{ if $i_{\rho \sigma} \equiv 1 \mod 2$, } \\
0  &\mod 2 &  \text{ otherwise. }
\end{array}
\right.
\]
These matrix elements define the map $\partial^{(0)}:C^{(0)}_i \rightarrow C^{(0)}_{i-1}$.

\begin{lem} \label{lem bd bd}
$\partial^{(0)} \circ \partial^{(0)} = 0$.
\end{lem}

For generators $[\rho] \in C_{i}^{(0)} = CF_{2i}$ and $[\tau] \in C^{(0)}_{i-2} = CF_{2i-4}$, we have
\[
< \bd \circ \bd ([\rho]), [\tau]> =
\sum_{[\sigma]} < \bd([\rho]), [\sigma]> <\bd([\sigma]), [\tau]>.
\]
Here $[\sigma]$ runs over the set of generators of $\cm_{i-1}$. 

If $i_{\rho \tau} \equiv 1 \mod 2$, then $i_{\rho \sigma}$ or $i_{\sigma \tau}$ is even by the additivity of the index. By definition, $<\bd ([\rho]), [\sigma] >$ or $<\bd ([\sigma]), [\tau] >$ is trivial, and hence $< \bd \circ \bd ([\rho]), [\tau] > \equiv 0 \mod 2$.

To prove the lemma in the case when $i_{\rho \tau} \equiv 0 \mod 2$, we consider the moduli space $M_{\rho \tau}'$. By the formula in Lemma \ref{lem 1}, 
\[
\dim M_{\rho \tau} 
\equiv \delta_{Y}([\rho]) - \delta_{Y}([\tau]) - 1 
\equiv 4 - 1 
\equiv 3 \mod 8.
\]
Hence we have the 2-dimensional moduli space $M_{\rho \tau}'$ of instantons with center of mass $0$. We need a real line bundle $\Lambda$ over the moduli space $M_{\rho \tau}'$, which is defined as in \cite{DK} for closed 4-manifolds. There is the universal bundle $\tilde{\E}$ over $(Y \times \R) \times \tilde{M}_{\rho \tau}$:
\[
\tilde{\E} := E \times_{\cG_0} \cA_{\rho \tau}^{\inst}
\longrightarrow
(Y \times \R) \times \tilde{M}_{\rho \tau}.
\]
Here $\cA_{\rho \tau}^{\inst}$ is the space of instantons with limits $\rho, \tau$ and $\cG_0$ is the group of gauge transformations with limit $1$ at $\pm \infty$.
For each $A \in \cA_{\rho \tau}^{\inst}$, we have the real part of the twisted Dirac operator $\Dir_{A}$:
\[
(\Dir_{A})_{\R}: 
L^{2,(-\alpha, \alpha)}_1( (S^+ \otimes E)_{\R} ) \longrightarrow L^{2, (-\alpha, \alpha)} ( (S^- \otimes E)_{\R} ).
\]
The universal bundle and the family of the real operators define the determinant line bundle over the framed moduli space:
\[
\tilde{ \Lambda } \stackrel{\R}{\longrightarrow} \tilde{M}_{\rho \tau}.
\]
We have a natural action of $\Gamma_{\rho} \times \Gamma_{\tau}$ on $\tilde{\Lambda}$ which is a lift of the action on $\tM_{\rho \tau}$. The subgroup $\{ \pm (1,1) \}$ acts on $\tM_{\rho \tau}$ trivially and on the fiber of $\tilde{\Lambda}$ with weight $\ind \Dir_A = i_{\rho \tau}$. Since we assumed that $i_{\rho \tau}$ is even, the action on the fiber is also trivial. Hence $\tilde{\Lambda}$ descends to a line bundle over $M_{\rho \tau}$. We denote the restriction to $M_{\rho \tau}'$ by $\Lambda$.

Take a generic section $s$ of $\Lambda$. We consider the end of the zero locus $s^{-1}(0)$. A standard argument, which can be found in \cite{D Floer}, shows the following:

\begin{lem}
Any sequence $\{ [A_{\alpha}] \}_{\alpha}$ in $M_{\rho \tau}'$ has a subsequence $\{ [A_{\alpha'}] \}_{\alpha'}$ which is chain-convergent to some $([A_1], [A_2]) \in M_{\rho \sigma}' \times M_{\sigma \tau}'$. Here $\sigma$ is a flat connection with $\Gamma_{\sigma} \cong U(1)$ and $\dim M_{\rho \sigma}' = \dim M_{\sigma \tau}' = 0$.
\end{lem}

As is well known, gluing instantons gives the map
\begin{equation*}
Gl:\coprod_{[\sigma]} M_{\rho \sigma}' \times M_{\sigma \tau}' \times 
\big( \Gamma_{\sigma} / \{ \pm 1 \} \big) \times (T_0, \infty) 
\longrightarrow
M_{\rho \tau}',
\end{equation*}
where $[\sigma]$ runs over the set of generators of $\cm_{i-1}$, and $T_0$ is a large positive number. 
The gluing map $Gl$ is a homeomorphism to an open set in $M_{\rho \tau}'$ and the complement of the image of $Gl$ is compact. 

Fix $T_1 > T_0$ and put $M_{\rho \tau}'' := M_{\rho \tau}' \backslash \im Gl_{>T_1}$. Here $Gl_{>T_1}$ is the restriction of $Gl$ to the domain where the parameter $T$ is larger than $T_1$. For a generic section $s$ of $\Lambda$, $N_{\rho \tau}'' := s^{-1}(0) \cap M_{\rho \tau}''$ is a smooth compact 1-dimensional manifold with boundary
\[
\coprod_{[\sigma]} \big( M_{\rho \sigma}' \times M_{\sigma \tau}' \times \big( \Gamma_{\sigma} / \{ \pm 1 \} \big) \times \{ T_1 \} \big) \cap s^{-1}(0).
\]
For $\fa \in M_{\rho \sigma}' \times M_{\sigma \tau}'$, we denote by $U(1)_{\fa}$ the corresponding gluing parameter. That is, $U(1)_{\fa} := \{ \fa \} \times \big( \Gamma_{\sigma} / \{ \pm 1 \} \big) \times \{ T_1 \} \cong U(1)$.

\begin{lem} \label{lem Lambda}
The line bundle $\Lambda$ is non-trivial on $U(1)_{\fa}$ if and only if 
\[
i_{\rho \sigma} \equiv  1 \mod 2.
\]
\end{lem}

Note that $i_{\rho \sigma} \equiv i_{\sigma \tau} \mod 2$ since we assumed $i_{\rho \tau} \equiv 0 \mod 2$.

This lemma can be proved in the same way as \cite[Lemma 3.14]{S Floer}. 
We give outline of the proof. Let $p:\Gamma_{\sigma} \rightarrow U(1)_{\fa}$ be the projection. The gluing theory gives a natural trivialization
\[
\hat{\Lambda} := p^* \Lambda|_{U(1)_{\fa}} \cong \underline{\R}.
\]
We have the natural action of $\Z_2$ on $\hat{\Lambda}$ and $\hat{\Lambda}/\Z_2 = \Lambda|_{U(1)_{\fa}}$. Through the trivialization, the action of $-1 \in \Z_2$ on the fiber is $(-1)^{i_{\rho \sigma}}$. ( See  (\ref{eq g -g}) below.) Hence we have obtained the statement.

\vspace{2mm}

By this lemma, we have
\[
\begin{split}
\# \partial N_{\rho \tau}''
&\equiv \sum_{[\sigma]; \ i_{\rho \sigma} \equiv i_{\sigma \tau} \equiv 1 \mod 2} 
\# M_{\rho \sigma}' \cdot \# M_{\sigma \tau}' \mod 2 \\
&\equiv \sum_{[\sigma]} < \bd ([\rho]), [\sigma]> < \bd([\sigma]), [\tau]> \mod 2 \\
&\equiv < \bd \circ \bd ([\rho]), [\tau]> \mod 2.
\end{split}
\]
On the the hand, the number of the boundaries of a 1-dimensional compact manifold is even. Hence we have obtained the required identity
\[
<\bd \circ \bd ([\rho]), [\tau]> \equiv 0 \mod 2.
\]
We can also show $\partial^{(1)} \circ \partial^{(1)} = 0$ by the same arguments.

\begin{dfn}
$I^{(0)}(Y;\fs):=H_*(\cm_*(Y), \bd)$, $I^{(1)}(Y;\fs) := H_*(C_*^{(1)}(Y), \partial^{(1)})$.
\end{dfn}

\subsection{Analog of Fukaya-Floer homology} \label{ss FF}

The aim of this subsection is to construct analog of Fukaya-Floer homology \cite{Fuk} for lens spaces. Fukaya-Floer homology is defined for a triple of a 3-manifold $Z$, a $U(2)$-bundle $Q$ over $Z$ and a loop $\gamma$ in $Z$, provided that all (projectively) flat connections on $Q$ are irreducible. The boundary operator is defined using not only 0-dimensional moduli spaces over $Z \times \R$ but also 2-dimensional moduli spaces. As in the previous case, to extend this construction to lens spaces $Y$, we must change the definition of the boundary using twisted Dirac operators over $Y \times \R$. Furthermore, as we will see later, we need to look out the contribution of the trivial connection differently from \cite{Fuk}, \cite{BD} and the previous subsection. The discussion which involves the trivial flat connection is similar to that in \cite{S Floer}.

\vspace{3mm}

Throughout this subsection, we assume $p$ is an odd prime integer. The assumption that $p$ is odd implies that $Y$ has an unique spin structure (up to isomorphism), and we have
\[
\delta_{Y}([\rho]) \equiv 0 \mod 2
\]
for all flat connections. See Corollary \ref{cor even}. Moreover the only trivial flat connection has $SU(2)$ as the stabilizer in the gauge group. 

\vspace{2mm}

Let $\gamma$ be a simple closed curve in $Y$. Put 
\[
C_{i}(Y;\gamma):= 
\left\{
\begin{array}{cl}
CF_{2i}(Y) \oplus CF_{2i-2}(Y) & \text{ if $i \not \equiv 0 \mod 8$, } \\
CF_{0}(Y) \oplus CF_{-2}(Y) \oplus \Z_2< [\theta] > & \text{ if $i \equiv 0 \mod 8$ }.
\end{array}
\right.
\]
We will define the boundary operator $\bdg:C_{*}(Y;\gamma) \rightarrow C_{*-1}(Y;\gamma)$ as follows.

As before, we define the matrix elements
\[
< \bdg([\rho]), [\sigma]> \in \Z_2
\]
for generators $[\rho] \in C_i(Y;\gamma), [\sigma] \in C_{i-1}(Y;\gamma)$ using moduli spaces over $Y \times \R$. First assume that 
\[
[\rho], [\sigma] \not = [\theta], \quad
\delta_Y([\rho]) - \delta_{Y}([\sigma]) \equiv 2 \mod 8.
\]
In this case, we have the moduli space $M_{\rho \sigma}'$ of dimension $0$. As before we define
\[
< \bdg([\rho]), [\sigma] > :=
\left\{
\begin{array}{cll}
\# M_{\rho \sigma}' &\mod 2 & \text{ if $i_{\rho \sigma} \equiv 1 \mod 2$}, \\
0 &\mod 2 & \text{otherwise}.
\end{array}
\right.
\]
Next we consider the case where
\[
[\rho], [\sigma] \not= [\theta], \quad
\delta_{Y}([\rho]) - \delta_{Y}( [\sigma] ) \equiv 4 \mod 8.
\]
In this case we have the moduli space $M_{\rho \sigma}'$ of dimension $2$. To define the matrix element, we use the determinant line bundle of twisted $\bar{\partial}$-operators over $\gamma \times \R$. Take a spin structure of $\gamma$ which represent the trivial spin bordism class. ( See \cite[Remark 2.3]{S Floer} for the reason why we choose the spin structure. ) The spin structure induces a spin structure of $\gamma \times \R$ (i.e. a square root $K_{\gamma \times \R}^{ \frac{1}{2} }$ of the canonical line bundle $K_{\gamma \times \R}$ ) and we have twisted $\bar{\partial}$-operators 
\[
\bar{\partial}_A:
\Gamma(K_{ \gamma \times \R }^{ \frac{1}{2} } \otimes E) \longrightarrow
\Gamma(K_{ \gamma \times \R }^{ \frac{1}{2} } \otimes E \otimes \Lambda^{0,1}_{\gamma \times \R}).
\]
Here $E$ is the rank 2 complex vector bundle over $Y \times \R$ associated with $\pi^*Q$ and $A$ is a connection on $\pi^* Q$ with limit $\rho, \sigma$. As in the previous subsection, the family of twisted $\bar{\partial}$ operators defines the determinant line bundle $\tilde{\cL}_{\gamma; \rs}$ over the framed moduli space $\tM_{\rs}'$. We show that this line bundle descends to a line bundle over $M_{\rho \sigma}'$. It is sufficient to prove that the index of $\bar{\partial}_A$ is even.

\begin{lem} \label{lem ind even}
When $p$ is prime, we have:
\begin{enumerate}
\item \label{trivial}
Assume that $[\gamma] \in H_1(Y;\Z)$ is trivial. Then for any flat connections $\rho, \sigma$  connections $A$ with limits $\rho, \sigma$, the index $\ind^{-+} \bar{\partial}_{A}$ is even.

\item
Assume that $[\gamma] \in H_1(Y;\Z)$ is non-trivial. Take flat connections $\rho, \sigma$ with $\Gamma_{\rho}, \Gamma_{\sigma} \cong U(1)$.
For any connection $A$ over $Y \times \R$ with limits $\rho, \sigma$, the index $\ind^{-+} \bar{\partial}_{A}$ is even.
\end{enumerate}
\end{lem}

We can see (\ref{trivial}) in this lemma as follows. Since $[\gamma] \in H_1(Y;\Z)$ is trivial, the restrictions $\rho|_{\gamma}, \sigma|_{\gamma}$ are gauge equivalent to the trivial flat connection over $\gamma$. In particular, $\rho|_{\gamma}$ is gauge equivalent to $\sigma|_{\gamma}$. By additivity of the index, $\ind^{-+} \bar{\partial}_A$ is equal to the index of $\bar{\partial}$-operator over $\gamma \times S^1$ twisted by the $SU(2)$-bundle. This index is zero by the index theorem.

\vspace{2mm}

The second part of Lemma \ref{lem ind even} follows from:

\begin{lem} \label{lem ind odd}
Assume that $[\gamma] \in H_1(Y;\Z)$ is non-trivial and that $p$ is prime. Take flat connections $\rho, \sigma$ with $\Gamma_{\rho}, \Gamma_{\sigma} \cong U(1)$. For any connection $A$ with limits $\theta$ and $\rho$ or with limits $\sigma$ and $\theta$, the index $\ind^{-+} \bar{\partial}_{A}$ is odd.
\end{lem}

Assuming Lemma \ref{lem ind odd}, we give the proof of Lemma \ref{lem ind even}. By the additivity of the index, we have
\[
\ind^{-+} \bar{\partial}_{A_{ \theta \rho }} + \ind^{-+} \bar{\partial}_{ A_{ \rho \sigma } } = \ind^{-+} \bar{\partial}_{ A_{\theta \sigma} }.
\]
Here $A_{\theta \rho}$ is a connection with limits $\theta, \rho$. Similarly for $A_{\rho \sigma}, A_{\theta \sigma}$. By Lemma \ref{lem ind odd}, both $\ind^{-+} \bar{\partial}_{ A_{\theta \rho} }$ and $\ind^{-+} \bar{\partial}_{A_{ \theta \sigma }}$ are odd. Therefore $\ind^{-+} \bar{\partial}_{ A_{\rho \sigma }}$ is even.

We will give the proof of Lemma \ref{lem ind odd} at the end of this subsection.
\vspace{2mm}

We go back to the definition of $\bdg$. As before suppose that $\delta_{Y}([\rho]) - \delta_{Y}([\sigma]) \equiv 4 \mod 8$ and that $[\rho], [\sigma] \not= [\theta]$. By Lemma \ref{lem ind even}, we have the determinant line bundle
\[
\cL_{\gamma;\rs} := \tilde{\cL}_{\gamma; \rs} / \Gamma_{\rho} \times \Gamma_{\sigma}
\stackrel{\C}{\longrightarrow}
M_{\rs}'.
\]
We want to define the matrix element $< \bdg([\rho]), [\sigma] >$ to be $\# s_{\gamma}^{-1}(0) \mod 2$ for a generic section $s_{\gamma}:M_{\rs}' \rightarrow \cL_{\gamma; \rs}$ if $i_{\rho \sigma} \equiv 1 \mod 2$ and zero otherwise. However $M_{\rho \sigma}'$ is not compact in general. We need a section of $\cL_{\gamma;\rs}$ which is non-vanishing on the end of $M_{\rs}'$ and transverse to the zero section. The end can be described as in the proof of Lemma \ref{lem bd bd}. The end is the image of the gluing map
\[
Gl: \coprod_{ [\mu] } 
M_{\rho \mu}' \times M_{\mu \sigma}' \times 
( \Gamma_{ \mu }/ \{ \pm 1 \}) \times (T_0, \infty)
\longrightarrow
M_{\rho \sigma}'.
\]
Here $[\mu]$ runs over the the set of gauge equivalence of flat connections with 
$\delta_{Y}([\rho]) - \delta_{Y}([\mu]) \equiv 2 \mod 8$ and with $\Gamma_{\mu} \cong U(1)$. We can take a desired section $s_{\gamma}$ as follows. 

First for each $[\mu]$ fix generic sections $s_{\gamma: \rho \mu}$ and $s_{\gamma;\mu \sigma}$ of the determinant line bundles over $M_{\rho \mu}'$ and $M_{\mu \sigma}'$. The zero loci of these sections are empty.

Next we consider the end of $M_{\rho \sigma}'$ described by gluing instantons $[A_1] \in M_{\rho \mu}'$ and $[A_2] \in M_{\mu \sigma}'$, which is identified with $E= \big( \Gamma_{\mu} / \{ \pm 1 \} \big) \times (T_0, \infty)$. 
Let $\hat{\cL}$ be the pull-back of the restriction $\cL_{\gamma:\rho \sigma}|_{E}$ by the projection 
\[
\Gamma_{\mu} \times (T_0, \infty) \longrightarrow E= \big( \Gamma_{\mu} /\{ \pm 1 \} \big) \times (T_0, \infty).
\]
On $\Gamma_{\mu} \times (T_0, \infty)$, the additivity of the index gives an isomorphism
\begin{equation} \label{eq gl line}
\hat{Gl}:
(\cL_{\gamma:\rho \mu})_{[A_1]} \boxtimes (\cL_{\gamma: \mu \sigma})_{[A_2]}  
\stackrel{\cong}{\longrightarrow}
\hat{\cL}.
\end{equation}
For $(g, T) \in \Gamma_{\mu} \times (T_0, \infty)$, we have
\begin{equation} \label{eq g -g}
\hat{Gl}|_{(g,T)} = (-1)^{ \ind^{-+} \bar{\partial}_{A_1} } \cdot \hat{Gl}|_{(-g,T)}.
\end{equation}
Here $\hat{Gl}|_{(g,T)}$ is the restriction of $\hat{Gl}$ to the fiber over $(g, T)$.
This can be seen as follows. The gauge equivalence class of instanton corresponding to $(g, T)$ is obtained by gluing instantons $u(A_1)$ and $A_2$, where $u$ is a gauge transformation over $Y \times \R$ with limit $g$ at $+\infty$. On the other hand, the gauge equivalence class corresponding to $(-g,T)$ is obtained by gluing $-u(A_1)$ and $A_2$. The action of $-1$ on the fiber $(\cL_{\gamma; \rho \mu})_{[A_1]}$ is $(-1)^{\ind^{-+} \bar{\partial}_{A_1}}$. Hence we have obtained (\ref{eq g -g}).

By Lemma \ref{lem ind even}, $\ind^{-+} \bar{\partial}_{A_1} \equiv 0 \mod 2$. Thus we obtain:

\begin{lem} \label{lem gl line}
The above isomorphism (\ref{eq gl line}) descends to an isomorphism
\[
\cL_{\gamma: \rho \sigma}|_{E} \cong 
(\cL_{\gamma: \rho \mu})_{[A_1]} \boxtimes ( \cL_{\gamma:\mu \sigma} )_{[A_2]}
\]
over $E= \big( \Gamma_{\mu} /\{ \pm 1 \} \big) \times (T_0, \infty)$.
\end{lem}

As in \cite[Section 2]{S Floer}, we can construct a section $s_{\gamma}:M_{\rho \sigma}' \rightarrow \cL_{\rho \sigma}$ which is compatible with the identification of Lemma \ref{lem gl line}. That is, if $[A^{\alpha}]$ be a sequence of points in $M_{\rs}'$ converging to some $([A_1], [A_2]) \in M_{\rho \mu}' \times M_{\mu \sigma}'$ then $s_{\gamma}([A_{\alpha}]) \rightarrow s_{\gamma; \rho \mu}([A_1]) \otimes s_{\gamma;\mu \sigma}([A_2])$ in the sense of \cite[Definition 2.7]{S Floer}. The section $s_{\gamma}$ does not vanish on the ends of the moduli space, since $s_{\gamma;\rho \mu}$ and $s_{\gamma;\mu \sigma}$ are non-vanishing sections. Thus the zero locus $s_{\gamma}^{-1}(0)$ is compact. Perturbing $s_{\gamma}$ over a compact set in $M_{\rho \sigma}'$, we may assume that $s_{\gamma}$ is transverse to the zero section. Therefore $s_{\gamma}^{-1}(0)$ is a finite set. We define the matrix element by
\[
< \bdg([\rho]), [\sigma]> :=
\left\{
\begin{array}{cll}
\# s_{\gamma}^{-1}(0) & \mod 2 & \text{ if $i_{\rho \sigma} \equiv 1 \mod 2$ },\\
0 & \mod 2 & \text{ otherwise. }
\end{array}
\right.
\]

\vspace{2mm}

Next we define the terms which involve the trivial connection. Let $[\rho] \in CF_{2}(Y) \subset C_1(Y;\gamma)$, $[\sigma] = [\theta] \in C_0(Y;\gamma)$. Then we have a 0-dimensional moduli space $M_{\rho \theta}'$. We define
\[
< \bdg([\rho]), [\theta] > :=
\left\{
\begin{array}{cll}
\# M_{\rho \theta}' & \mod 2 & \text{ if $i_{\rho \theta} \equiv 1 \mod 2$, $[\gamma] \not= 0$ in $H_1(Y;\Z)$, } \\
0  & \mod 2 & \text{ otherwise. }
\end{array}
\right.
\]

\vspace{2mm}

Let $[\rho] = [\theta] \in C_0(Y;\gamma)$, $[\sigma] \in CF_{-4}(Y) \subset CF_{-1}(Y;\gamma)$. Then
\[
\dim M_{\theta \sigma}
\equiv \delta_{Y}([\theta]) - \delta_{Y}([\sigma]) - \dim \Gamma_{\theta}
\equiv 0 - (-4) - 3
\equiv 1 \mod 8.
\]
We have a $0$-dimensional moduli space $M_{\theta \sigma}'$. We put
\[
< \bdg([\theta]), [\sigma]> :=
\left\{
\begin{array}{cll}
\# M_{\theta \sigma}' & \mod 2 & \text{ if $i_{\theta \sigma} \equiv 1 \mod 2$, $[\gamma] \not= 0$ in $H_1(Y;\Z)$, } \\
0 & \mod 2 & \text{otherwise}.
\end{array}
\right.
\]
We define other matrix elements to be zero.

\vspace{2mm}

The part of the boundary map which does not involve the trivial flat connection is as the following diagram:
\[
\xymatrix@R=6pt{
C_i(Y;\gamma) \ar@{=}[d]  \ar[r]^{ \partial_{\gamma} } & C_{i-1}(Y;\gamma) \ar@{=}[d] \ar[r]^{ \partial_{\gamma} } & C_{i-2}(Y;\gamma) \ar@{=}[d] \\
CF_{2i}(Y)  \ar[r] \ar[rdd]  & CF_{2i-2}(Y)  \ar[r] \ar[rdd] & CF_{2i-4}(Y) \\
\oplus & \oplus & \oplus \\
CF_{2i-2}(Y) \ar[r] & CF_{2i-4}(Y) \ar[r] & CF_{2i-6}(Y)
}
\]
The horizontal maps are defined using 0-dimensional moduli spaces, and the diagonal maps are defined using 2-dimensional moduli spaces and the determinant line bundles of $\gamma \times \R$.

The part of the boundary map which involves the trivial flat connection is as the following diagram:
\begin{equation} \label{eq bd map}
\xymatrix@R=6pt{
C_{1}(Y;\gamma) \ar@{=}[d] \ar[r]^{ \partial_{\gamma} } & C_{0}(Y;\gamma) \ar@{=}[d] \ar[r]^{ \partial_{\gamma} } & C_{-1}(Y;\gamma) \ar@{=}[d] \\
CF_2(Y) \ar[r]^{a} \ar[rdd]^{c} \ar[rdddd]^{e} & CF_{0}(Y) \ar[r] \ar[rdd]^{b} & CF_{-2}(Y) \\
\oplus & \oplus & \oplus \\
CF_{0}(Y) \ar[r] & CF_{-2}(Y) \ar[r]^{d} & CF_{-4}(Y) \\ 
                 & \oplus            &            \\
                 & \Z_2 <[\theta] > \ar[ruu]_f & 
}
\end{equation}
The $a, c$ and $e$ are maps from $CF_2(Y)$, and the $b, d$ and $f$ are maps into $CF_{-4}(Y)$. (These notations will be used later.)

\begin{lem} \label{lem bdg bdg}
$\bdg \circ \bdg = 0$.
\end{lem}

We must show $< \bdg \circ \bdg([\rho]), [\tau]> \equiv 0 \mod 2$ for all generators $[\rho] \in C_i(Y;\gamma), [\tau] \in C_{i-2}(Y;\gamma)$. We give the proof in the case where $[\rho] \in CF_{2i}(Y) \subset C_{i}(Y;\gamma)$, $[\tau] \in CF_{2i-6}(Y) \subset C_{i-2}(Y;\gamma)$. The proof for the other cases is the same as that of Lemma \ref{lem bd bd}. 

If $i_{\rho \tau} \equiv 1 \mod 2$, then we have $i_{\rho \sigma} \equiv 1 \mod 2$ or $i_{\sigma \tau} \equiv 1 \mod 2$ for generators $[\sigma]$ of $C_{i-1}(Y;\gamma)$. Hence $< \bdg([\rho]), [\sigma] > \equiv 0 \mod 2$ or $< \bdg([\sigma]).[\tau] > \equiv 0 \mod 2$ by definition. Since 
\[
< \bdg([\rho]), [\tau]> = \sum_{[\sigma]} < \bdg([\rho]), [\sigma]> < \bdg ([\sigma]), [\tau]>,
\]
we have $< \bdg \circ \bdg ([\rho]), [\tau] > \equiv 0 \mod 2$.

Suppose that $i_{\rho \tau} \equiv 0 \mod 2$. We consider the moduli space $M_{\rho \tau}'$ of dimension $4$. We analyze the end of a 1-dimensional moduli space
\[
N=M_{\rho \tau}' \cap s_{\gamma; \rho \tau}^{-1}(0) \cap s_{\Lambda; \rho \tau}^{-1}(0),
\]
where $s_{\gamma; \rho \tau}$ and $s_{\Lambda; \rho \tau}$ are sections of $\cL_{\gamma; \rho \tau} \rightarrow M_{\rho \tau}'$ and $\Lambda \rightarrow M_{\rho \tau}'$ respectively. 

A dimension counting argument shows the following:

\begin{lem} \label{lem end 2}
Let $\{ [A^{\alpha}] \}_{\alpha}$ be a sequence in $M_{\rho \tau}'$. Then we can find a subsequence $\{ [A^{\alpha'}] \}_{\alpha'}$ such that
\begin{enumerate}[(i)]
\item
$[A^{\alpha'}] \longrightarrow ([A_1],[A_2]) \in M_{\rho \sigma}' \times M_{\sigma \tau}'$, where $\sigma$ is a flat connection with $\Gamma_{\sigma} \cong U(1)$,and $\dim M_{\rho \sigma}'=2$, $\dim M_{\sigma \tau}'=0$, (i.e. $[\sigma] \in CF_{2i-4}$) or

\vspace{2mm}

\item
$[A^{\alpha'}] \longrightarrow ([A_1],[A_2]) \in M_{\rho \sigma}' \times M_{\sigma \tau}'$, where $\sigma$ is a flat connection with $\Gamma_{\sigma} \cong U(1)$,and $\dim M_{\rho \sigma}'=0$, $\dim M_{\sigma \tau}'=2$, (i.e. $[\sigma] \in CF_{2i-2}$ ) or

\vspace{2mm}

\item \label{trivial flat}
$[A^{\alpha'}] \longrightarrow ([A_1], [A_2]) \in M_{\rho \theta}' \times M_{\theta \sigma}'$, and $\dim M_{\rho \theta}'= 0, \dim M_{\theta \sigma}'= 0$.

\end{enumerate}
\end{lem}

Note that (\ref{trivial flat}) occurs only if $[\rho] \in CF_{2}(Y) \subset CF_{1}(Y;\gamma)$, $[\tau] \in CF_{-4}(Y) \subset CF_{-1}(Y;\gamma)$, and that the case where $\dim M_{\rho \sigma}' = \dim M_{\sigma \tau}' = 1$ does not occur since $\delta_{Y}([\rho])-\delta_{Y}([\sigma]) \equiv 0 \mod 2$ for all $\rho, \sigma$ as we will prove in subsection \ref{ss moduli space}. (Corollary \ref{cor even})

\vspace{2mm}

We consider the case where (\ref{trivial flat}) may occur. That is, $[\rho] \in CF_{2}(Y)$, $[\tau] \in CF_{-4}(Y)$. As usual, we take a section $s_{\gamma: \rho \tau}$ of $\cL_{\gamma;\rho \tau} \rightarrow M_{\rho \tau}'$ such that if $[A^{\alpha}] \rightarrow ([A_1], [A_2]) \in M_{\rho \sigma}' \times M_{\sigma \tau}'$ and $\Gamma_{\sigma} \cong U(1)$, then $s_{\gamma;\rho \tau}([A^{\alpha}]) \rightarrow s_{\gamma;\rho \sigma}([A_1]) \otimes s_{\gamma;\sigma \tau}([A_2])$. By Lemma \ref{lem end 2}, the end of $N$ is identified with
\begin{equation} \label{eq end N 2}
\begin{split}
& \coprod_{
\begin{subarray}{c}
[\sigma]: \\
\delta_{Y}([\sigma]) \equiv -2 \mod 8
\end{subarray}
}
\coprod_{\fa}
\left(
U(1)_{\fa} \cap s^{-1}_{\Lambda; \rho \tau}(0)
\right) \cup 
\\
& \coprod_{
\begin{subarray}{c}
[\sigma]: \\
\delta_{Y}([\sigma]) \equiv 0 \mod 8 \\
\Gamma_{\sigma} \cong U(1)
\end{subarray}
}
\coprod_{\fb}
\left(
U(1)_{\fb} \cap  s^{-1}_{\Lambda: \rho \tau}(0)
\right) \cup 
\\
&\coprod_{\fc}
\left(SO(3)_{\fc} \cap s_{\gamma: \rho \tau}^{-1}(0) \cap s_{\Lambda; \rho \tau}^{-1}(0)
\right).
\end{split}
\end{equation}
Here $\fa$, $\fb$ and $\fc$ run over
\[
\big( M_{\rho \sigma}' \cap s_{\gamma;\rho \sigma}^{-1}(0) \big) \times M_{\sigma \tau}', \ 
M_{\rho \sigma}' \times \big( M_{\sigma \tau}' \cap s_{\gamma; \sigma \tau}^{-1}(0) \big), 
\ \text{and} \ 
M_{\rho \theta}' \times M_{\theta \tau}'
\]
respectively. For $\fa=([A_1], [A_2])$, $U(1)_{\fa}$ is the gluing parameter used to glue instantons $[A_1]$ and $[A_2]$. Similarly for $U(1)_{\fb}$ and  $SO(3)_{\fc}$. We denote the first term in (\ref{eq end N 2}) by $\partial_1 N$, and the second and third terms by $\partial_2 N$ and $\partial_3 N$. From Lemma \ref{lem Lambda} we have
\[
\# \big( U(1)_{\fa} \cap s_{\Lambda;\rho \tau}^{-1}(0) \big) =
\left\{
\begin{array}{cll}
1 & \mod 2 & \text{ if $i_{\rho \sigma} \equiv i_{\sigma \tau} \equiv 1 \mod 2$ } \\
0 & \mod 2 & \text{ otherwise}.
\end{array} 
\right.
\]
Therefore we obtain
\[
\begin{split}
\# \partial_1 N 
&\equiv
\sum_{
\begin{subarray}{c}
[\sigma]: \\
\delta_{Y}([\sigma]) \equiv -2 \mod 8 \\
i_{\rho \sigma} \equiv i_{\sigma \tau} \equiv 1 \mod 2
\end{subarray}
}
\# (M_{\rho \sigma}' \cap s_{\gamma;\rho \sigma}^{-1}(0)) \times M_{\sigma \tau}' \\
&\equiv 
\sum_{
\begin{subarray}{c}
[\sigma]: \\
\delta_{Y}([\sigma]) \equiv -2 \mod 8
\end{subarray}
}
< \bdg([\rho]), [\sigma]> < \bdg([\sigma]), [\tau]> \\
&\equiv
\sum_{[\sigma]}< a([\rho]), [\sigma]><b([\sigma]), [\tau]>
\mod 2.
\end{split}
\]
Here $a, b$ are the maps in the diagram (\ref{eq bd map}).
Similarly we have
\[
\begin{split}
\# \partial_2 N 
&\equiv
\sum_{
\begin{subarray}{c}
[\sigma]: \\
\delta_{Y}([\sigma]) \equiv 0 \mod 8 \\
\Gamma_{\sigma} \cong U(1)
\end{subarray}
}
< \bdg([\rho]), [\sigma]> < \bdg([\sigma]), [\tau]> \\
&\equiv
\sum_{[\sigma]}
< c([\rho]), [\sigma]> < d([\sigma]), [\tau]>
\mod 2.
\end{split}
\]
Here $c, d$ are also the maps in (\ref{eq bd map}).

To compute $\# \partial_3 N_3 \mod 2$, we must know whether $\Lambda|_{ SO(3)_{\fc} }$ and $\cL_{\rho \tau}|_{SO(3)_{\fc}}$ are trivial or not. As \cite[Lemma 3.14]{S Floer}, we can show that $\Lambda|_{SO(3)_{\fc}}$ is non-trivial if and only if 
\[
i_{\rho \theta} \equiv i_{\theta \tau} \equiv 1 \mod 2.
\]
Similarly $\cL_{\gamma; \rho \tau}|_{SO(3)_{\fc}}$ is non-trivial if and only if
\[ 
\ind^{-+} \bar{\partial}_{A_1} \equiv \ind^{-+} \bar{\partial}_{A_2} \equiv 1 \mod 2.
\]
Here $\fc=([A_1], [A_2]) \in M_{\rho \theta}' \times M_{\theta \tau}'$. By Lemma \ref{lem ind even} and \ref{lem ind odd}, we obtain:

\begin{lem} \label{lem L SO3}
The line bundle $\cL_{\gamma:\rho \tau}$ is non-trivial over $SO(3)_{\fc}$ if and only if $[\gamma] \in H_1(Y;\Z)$ is non-trivial. 
\end{lem}

Therefore it follows that
\[
\begin{split}
\# \partial_3 N 
&\equiv \left\{
\begin{array}{cll}
\# M_{\rho \theta}' \cdot \# M_{\theta \tau}' & \mod 2 
& \text{ if $i_{\rho \theta} \equiv i_{\theta \tau} \equiv 1 \mod 2$,  
$[\gamma] \not= 0 \in H_1(Y;\Z)$, } \\
0 & \mod 2 & \text{ otherwise}.
\end{array}
\right. \\
&\equiv < e([\rho]), [\theta] > < f([\theta]), [\tau] > \mod 2.
\end{split}
\]
Here $e$ and $f$ are the maps in (\ref{eq bd map}).
Thus we have
\[
\begin{split}
&\# \partial_1 N + \# \partial_2 N + \# \partial_3 N \\
& \quad \equiv 
\sum_{[\sigma]}< a([\rho]), [\sigma]> < b([\sigma]), [\tau]> +
\sum_{[\sigma]}< c([\rho]), [\sigma]> < d([\sigma]), [\tau]>  \\
& \qquad \qquad + < e([\rho]), [\theta] > < f([\theta]), [\tau] > \\
& \quad \equiv
< \bdg \circ \bdg([\rho]), [\tau] > \mod 2.
\end{split}
\]
Since the number of the ends of a 1-dimensional manifold is even, the left hand side is even. Thus we have obtained $< \bdg \circ \bdg([\rho]), [\tau] > \equiv 0 \mod 2$ as required.

\vspace{2mm}

The proof for the other cases is similar and we omit the proof.

\begin{dfn}
Let $p$ be an odd, prime integer and $Y$ be $L(p,q)$ or $-L(p,q)$.
Define $I_*(Y;\gamma) := H_*( C_*(Y;\gamma), \bdg)$.
\end{dfn}

We can show that $I_*(Y;\gamma)$ is independent of the choice of sections of $\cL_{\gamma; \rho \sigma}$, up to canonical isomorphism, using standard arguments. 

\vspace{2mm}

It remains to prove Lemma \ref{lem ind odd}. 

\noindent
{\it Proof of Lemma \ref{lem ind odd}}

Let $\rho$ be a flat connection over $Y$ with $\Gamma_{\rho} \cong U(1)$ and take a connection $A$ over $Y \times \R$ with limits $\theta, \rho$. 
We can take a gauge transformation $g$ over $\gamma$ such that 
\[
g^*(\rho) = \rho_l \oplus -\rho_l
\]
over $\gamma$. Here $l$ is a positive integer with $1 \leq l \leq p-1$,  $\rho_l = \frac{2\pi l \sqrt{-1}}{p} d \varphi$, and $\varphi$ is a coordinate of $\gamma$ such that the restriction of the Riemannian metric of $Y$ to $\gamma$ is written as $d\varphi^{\otimes 2}$. Note that the restriction $\rho|_{\gamma}$ is not gauge equivalent to the trivial connection because $[\gamma] \not= 0$ in $H_1(Y;\Z) = \Z_p$ and we assumed that $p$ is prime. Hence $l$ is not zero. Since $SU(2)$ is simply connected, we can take a gauge transformation $\tilde{g}$ over $\gamma \times \R$ such that
\[
\tilde{g}|_{ \gamma \times (-\infty, -1) } = 1, \quad
\tilde{g}|_{ \gamma \times (1, \infty)} = g.
\]
Since $\ind^{-+} \bar{\partial}_{A}$ depends only on the limits of the restriction $A|_{\gamma \times \R}$, it is sufficient to consider a connection $A$ of the form
\[
A = a \oplus -a.
\]
Here $a$ is a ${\rm U(1)}$-connection such that
\[
a = 
\left\{
\begin{array}{cl}
\theta & \text{ on $\gamma \times (-\infty, -1)$ }, \\
\rho_{l} & \text{ on $\gamma \times (1, \infty)$ }.
\end{array}
\right.
\]
The index $\ind^{-+} \bar{\partial}_A$ is the sum $\ind^{-+} \bar{\partial}_{a} + \ind^{-+} \bar{\partial}_{-a}$. We compute $\ind^{-+} \bar{\partial}_{a}$ and $\ind^{-+} \bar{\partial}_{-a}$.

For $t \in [0, 1]$, put
\[
a_t = \frac{2 \pi l \sqrt{-1} t}{p} d\varphi.
\]
We give the complex structure $\gamma \times \R$ using the coordinate
\[
z = t + \sqrt{-1} \varphi
\]
where $t$ is the coordinate of $\R$ and $\varphi$ is the coordinate of $\gamma$.
We trivialize the line bundles $K_{\gamma \times \R} = \Lambda^{1,0}_{\gamma \times \R}$ and $\Lambda^{0,1}_{\gamma \times \R}$ using $e^{ \sqrt{-1} \varphi } dz = e^{ \sqrt{-1} \varphi } (dt + \sqrt{-1} d \varphi)$ and $d\bar{z} = dt - \sqrt{-1} d \varphi$ respectively. (The factor $e^{ \sqrt{-1} \varphi}$ makes the trivialization of $K_{\gamma \times \R}$ be compatible with the spin structure of $\gamma \times \R$ chosen before.)
Through these trivializations, the twisted $\bar{\partial}$-operator $\bar{\partial}_{a_t}$ is written as
\[
\bar{\partial}_{a_t} = 
\frac{1}{2} \left(
\frac{ \partial }{ \partial t} + \sqrt{-1} \frac{\partial}{\partial \varphi} -  \frac{2 \pi l  t}{p}
\right).
\]
The index $\ind^{-+} \bar{\partial}_{a}$ is equal to the spectral flow of the family
\[
\sqrt{-1} \frac{\partial}{\partial \varphi}  -  \frac{2 \pi l t}{p} + \epsilon
\quad
(0 \leq t \leq 1),
\]
where $\epsilon > 0$ is the small number used to define the weighted Sobolev spaces $L^{2, (-\epsilon, \epsilon)}$, $L^{2, (-\epsilon, \epsilon)}_1$. The spectra of this family are
\[
\lambda_n(t) = -\frac{2 \pi l t}{p} + \epsilon + 2 \pi n \quad (n \in \Z).
\]
>From this, we have $\ind^{-+} \bar{\partial}_{a} = -1$.

\vspace{2mm}

Similarly, the index $\ind^{-+} \bar{\partial}_{-a}$ is equal to the spectral flow of
\[
\sqrt{-1} \frac{\partial}{\partial \varphi} + \frac{2 \pi l t}{p} + \epsilon
\quad
( 0 \leq t \leq 1).
\]
It is easy to see that $\ind^{-+} \dbar_{-a_t} = 0$. Therefore we have $\ind^{-+}\bar{\partial}_{A} = -1$.

\vspace{2mm}

The proof for connections $A$ with limits $\sigma, \theta$ is similar.

\section{Gluing formula}

\subsection{2-torsion invariant for closed 4-manifolds} \label{ss torsion inv}
We show a gluing formula for an invariant $\Psi^{u_1}_X$ of non-spin closed 4-manifolds $X$ introduced in \cite{S torsion}. We recall the definition of $\Psi^{u_1}_X$ briefly.¡¡(See also \cite{FS} for the case when $X$ is spin.)
  This invariant is defined to be a function on a subspace of $\oplus_{ d \geq 0 } H_2(X;\Z)^{ \otimes d}$ as follows.

Let $X$ be a closed, non-spin, simply connected 4-manifold with $b^+ > 1$ and even. Take a principal $SO(3)$-bundle $P$ over $X$ with $w_2(P) = w_2(X)$ and with $p_1(P) \equiv \sigma(X) \mod 8$. Here $\sigma(X)$ is the signature of $X$. Fix a Riemannian metric on $X$. Then we have the moduli space $M_P=M_P(g)$ of instantons on $P$. For generic $g$, $M_P$ is a smooth manifold of dimension 
\[
-2p_1(P) - 3(1 + b^+(X)).
\]
Since we assumed $b^+(X)$ is even, the dimension is odd, and we can write $\dim M_{P}=2d+1$ for some integer $d$. 

Suppose $d \geq 0$ and take $d$ homology classes $[\Sigma_1], \dots, [\Sigma_d] \in H_2(X;\Z)$ with self-intersection number even. Then we have the determinant line bundles $\cL_{\Sigma_i} \rightarrow M_P$ of the twisted $\bar{\partial}$-operator over $\Sigma_i$. We can take sections $s_{\Sigma_i}$ of $\cL_{\Sigma_i}$ such that
\[
N=M_P \cap V_{\Sigma_1} \cap \cdots \cap V_{\Sigma_d}
\]
is a compact smooth manifold of dimension $1$. (See \cite{S torsion}.) Here $V_{\Sigma_i}$ is the zero locus of $s_{\Sigma_i}$.

Take a $U(2)$-bundle $Q$ over $X$ with $Q / U(1) = P$ and a spin-c structure $\fs_X$ of $X$ with $c_1(\det \fs ) = - c_1(Q)$. Then we have the real part
\[
( \Dir_A )_{ \R }:
\Gamma ( ( S^+ \otimes E )_{\R}) \longrightarrow \Gamma ( (S^- \otimes E)_{\R})
\]
of the twisted Dirac operator $\Dir_A$. Here $A$ is a connection on $Q$, $E$ is the rank 2 complex vector bundle associated with $Q$, and $S^{ \pm }$ is the spinor bundles of $\fs$. We denote by $\Lambda$ the determinant line bundle over $M_P$ of the family $\{ ( \Dir_A )_{\R} \}_{[A]}$.
We define $\Psi_{X}^{u_1}([\Sigma_1], \dots, [\Sigma_d]) \in \Z_2$ to be 
\[
\Psi_{X}^{u_1}([\Sigma_1], \dots, [\Sigma_d]) = \# N \cap s_{ \Lambda }^{-1}(0) \mod 2
\]
for a generic section $s_{\Lambda}$ of $\Lambda$. We can see that this is independent of the choices of the Riemannian metric and the sections of the line bundles.

\subsection{Relative invariants} \label{ss rel inv}

In this subsection, we generalize $\Psi_X^{u_1}$ to compact manifolds with boundary $Y$. Here $Y$ denote  $L(p,q)$ or $-L(p,q)$ as usual. Throughout this subsection we assume that $p$ is odd. 

\vspace{2mm}

Let $X_1$ be a compact, connected, simply connected, non-spin 4-manifold with boundary $Y$. Assume that $b^+(X_1) > 0$. We take a Riemannian metric on $X_1$ whose restriction to $Y$ coincides with the standard metric. Before we define the relative invariant of $X_1$, we discuss the dimension of moduli spaces of instantons over $\hat{X}_1 = X_1 \cup (Y \times \R_{\geq 0})$. To calculate the dimension, we need the following:

\begin{lem} \label{lem alpha}
Let $P_1$ be an $SO(3)$-bundle over $X_1$. There exists a cohomology class $\alpha \in H^2(X_1;\Z)$ with the following properties:

\begin{enumerate}[(i)]
\item
$\alpha \equiv w_2(P_1) \mod 2$,

\item
$\alpha|_Y = 0$ in $H^2(Y;\Z)$.

\end{enumerate}
\end{lem}

We will give the proof later. Fix a cohomology class $\alpha  \in H^2(X_1;\Z)$ with the above properties. From the exact sequence
\[
H^1(Y;\Z) = 0 \rightarrow H^2(X_1,Y;\Z) \rightarrow H^2(X_1;\Z) \rightarrow H^2(Y;\Z) \rightarrow \cdots
\]
and the property that $\alpha|_Y = 0$, there is the unique lift $\tilde{\alpha} \in H^2(X_1,Y;\Z)$ of $\alpha$. We define 
\[
\alpha^2 := < \alpha \cup \tilde{\alpha}, [X_1,Y]> \in \Z.
\]

\begin{prop} \label{prop moduli X_1 dim}
We denote by $\hat{P}_1$ the extension of $P_1$ to $\hat{X}_1 = X_1 \cup (Y \times \R_{ \geq 0 })$ . For any flat connection $\rho$ over $Y$, we have
\[
\dim M_{ \hat{P}_1, \rho} \equiv -\delta_Y([\rho]) - 2 \alpha^2 - 3(1+b^+(X_1)) \mod 8.
\]
Here $M_{ \hat{P}_1, \rho}$ is the moduli space of instantons on $\hat{P}_1$ with limit $\rho$.
\end{prop}

We prove Lemma \ref{lem alpha} and Proposition \ref{prop moduli X_1 dim}.
Lemma \ref{lem alpha} follows from:

\begin{lem} \label{lem sur}
The maps $H^2(X_1;\Z) \rightarrow H^2(X_1;\Z_2)$ and $H^2(X_1;\Z) \rightarrow H^2(Y; \Z)$ are surjective.
\end{lem}

Assuming this lemma, we prove Lemma \ref{lem alpha}. By Lemma \ref{lem sur}, we can find a lift $\alpha' \in H_2(X_1;\Z)$ of $w_2(P_1) \in H^2(X_1;\Z_2)$. By Bockstein exact sequence
\[
H^2(Y;\Z) \stackrel{\times 2}{\longrightarrow} H^2(Y;\Z) \rightarrow H^2(Y;\Z_2) = 0
\] 
there is an element $\beta' \in H^2(Y;\Z)$ such that $\alpha'|_Y = 2\beta'$. (Recall that we assumed that $p$ is odd.) By Lemma \ref{lem sur}, we have an extension $\alpha'' \in H^2(X_1;\Z)$ of $\beta' \in H^2(Y;\Z)$. Putting $\alpha := \alpha' - 2\alpha"$, we obtain a cohomology class with the required properties.

\vspace{2mm}

To prove Lemma \ref{lem sur}, we consider the exact sequences:
\begin{gather*}
\cdots \rightarrow H_1(Y;\Z) \rightarrow H_1(X_1;\Z) \rightarrow H_1(X_1,Y;\Z) \rightarrow 0  \\
\cdots \rightarrow H^1(X_1;\Z_2) \rightarrow H^2(X_1;\Z) \stackrel{\times 2}{\rightarrow} H^2(X_1;\Z) \rightarrow H^2(X_1;\Z_2) \rightarrow H^3(X_1;\Z) \rightarrow \cdots 
\end{gather*} 
>From the first sequence and the fact that $H_1(X_1;\Z) = 0$ we see that $H_1(X_1,Y;\Z) = 0$. (We assumed that $X_1$ is simply connected.) By Poincare duality, we also have $H^3(X_1;\Z) = 0$. From the second sequence, it follows that $H^2(X_1;\Z) \rightarrow H^2(X_1;\Z_2)$ is surjective.

We also see that $H^2(X_1;\Z) \rightarrow H^2(Y;\Z)$ is surjective from the exact sequence
\[
H^2(X_1,Y;\Z) \rightarrow H^2(X_1;\Z) \rightarrow H^2(Y;\Z) \rightarrow H^3(X_1,Y;\Z) \cong H_1(X_1;\Z) = 0.
\]

\vspace{2mm}

\noindent
{\it Proof of Proposition \ref{prop moduli X_1 dim}}

Choose a cohomology class $\alpha \in H^2(X_1;\Z)$ with the properties in Lemma \ref{lem alpha}.
Let $Q_1 \rightarrow X_1$ be a $U(2)$-bundle with $c_1(Q_1) = \alpha$ and fix an identification $Q_1/U(1) = P_1$. Since $\alpha|_{Y} = 0$, we can take a trivialization $\varphi$ of $Q_1|_Y$. We write $\eta$ for the trivial connection with respect to $\varphi$. The fact that $H^1(X_1;\Z_2) = 0$ implies that the moduli space $M_{ \hat{Q}_1,\eta }$ of instantons on $\hat{Q}_1$ with fixed determinant is naturally identified with $M_{ \hat{P}_1, \theta }$. Here $\theta$ is the trivial flat connection on $P_1|_{Y}$ with respect to the trivialization induced by $\varphi$.

Take a compact oriented 4-manifold $X_2$ with boundary $-Y$. Then we have a closed 4-manifold $X = X_1 \cup_{Y} X_2$. Using the trivialization $\varphi$, we extend $Q_1$ to $X$ in the obvious way. We write $Q_X$ for the $U(2)$-bundle over $X$. By the index formula, we have
\begin{equation} \label{eq QX}
\dim M_{Q_X} = 8c_2(Q_X) - 2 \alpha^2 - 3(1 - b_1(X) + b^+ (X)).
\end{equation}
Here we used the fact that $c_1(Q_X)^2 = \alpha^2$.  By the additivity of the index, we can also write
\begin{equation} \label{eq add}
\dim M_{Q_X} = \dim M_{\hat{Q}_1, \eta} + \dim M_{\hat{Q}_2, \eta} + 3.
\end{equation}
Here $\hat{Q}_2$ is the trivial $U(2)$-bundle over $\hat{X}_2$. The formal dimension $\dim M_{ \hat{Q}_2, \eta}$ is equal to the index of
\[
d^* + d^+:
L^{2,\epsilon}_1( \Lambda_{ \hat{X}_2}^1 \otimes \su(2))
\longrightarrow
L^{2, \epsilon}( (\Lambda_{ \hat{X}_2}^0 \oplus \Lambda_{ \hat{X}_2 }^+) \otimes \su(2)).
\]
It is easy to see that the index is
\begin{equation} \label{eq X2}
-3 (1 - b_1(X_2)  + b^+(X_2)).
\end{equation}
>From (\ref{eq QX}), (\ref{eq add}) and (\ref{eq X2}), we obtain 
\[
\begin{split}
\dim M_{ \hat{P}_1, \theta } 
&= \dim M_{ \hat{Q}_1, \eta} \\
&=8c_2(Q_X) - 2\alpha^2 - 3( 1 + b^+(X_1)) \\
&\equiv -2 \alpha^2 - 3(1 + b^+(X_1)) \mod 8.
\end{split}
\]
It follows from the additivity of the index that
\[
\dim M_{ \hat{P}_1, \rho} 
\equiv \dim M_{ \hat{P}_1, \theta} - \delta_{Y}([\rho])
\equiv -\delta_{Y}([\rho]) -2 \alpha^2 - 3(1 + b^+(X_1)) \mod 8
\]
for any flat connection $\rho$. 

\vspace{3mm}

From now on, we suppose that $w_2(P_1) \equiv w_2(X_1) \mod 2$, and fix a cohomology class $\alpha \in H^2(X_1;\Z)$ with the properties in Lemma \ref{lem alpha}. Take a $U(2)$-bundle $Q_1$ on $X_1$ as in the proof of Proposition \ref{prop moduli X_1 dim}. 
As before, we identify the moduli spaces of instantons on $\hat{P}_1$ with the moduli spaces of instantons on $\hat{Q}_1$ with fixed determinant $a_{\det}$. Here $a_{\det}$ is a fixed connection on the $U(1)$-bundle over $\hat{X}_1$ induced by $\hat{Q}_1$. Suppose that the limit of $a_{\det}$ at $\infty$ is the trivial connection.

Since we supposed $p$ is odd, $H^1(Y;\Z_2) = 0$. Hence the moduli space of $SO(3)$-flat connections over $Y$ and the moduli spaces of $SO(3)$-instantons over $Y \times \R$ are also naturally identified with the moduli spaces of the $SU(2)$-bundle. Hence we can regard Floer homologies in the previous section as those defined by using the $SO(3)$-bundle $P_1|_Y$. 

\vspace{2mm}

Suppose that $b^+(X_1)$ is odd. Then $-2\alpha^2 - 3 ( 1 + b^+(X_1))$ is even. Put
\[
d:= -\frac{ 2 \alpha^2 + 3 ( 1 + b^+ (X_1)) }{2}.
\]
We define relative invariants of $X_1$ in this situation. We consider three cases. In the first case, we will define $\Psi_{X_1}^{u_1} \in I_{d}(Y)$ using 0-dimensional moduli spaces over $\hat{X}_1$. (Since $\delta_Y \equiv 0 \mod 2$, $I_*^{(1)}(Y) = 0$. Hence we write $I_*(Y)$ for $I_*^{(0)}(Y)$.) In the second case, we consider a homology class $[\Sigma_1] \in H_2(X_1;\Z)$ with $[\Sigma_1] \cdot [\Sigma_1] \equiv 0 \mod 2$. Here $\Sigma_1$ is a closed oriented surface embedded in $X_1$. Using 2-dimensional moduli spaces over $\hat{X}_1$, we will define $\Psi^{u_1}_{X_1}([\Sigma_1]) \in I_{d-1}(Y)$. In the last case, we consider a relative homology class $[\Sigma_1] \in H_2(X_1,Y;\Z)$. Here $\Sigma_1$ is an embedded surface $X_1$ with $\partial \Sigma_1 = \gamma \subset Y$. We will define a relative invariant $\Psi_{X_1}^{u_1}([\Sigma_1]) \in I_{d}(Y,\gamma)$ using 0-dimensional and 2-dimensional moduli spaces, provided the class $[\Sigma_1]$ satisfies some conditions.

\vspace{3mm}

As in the definition of the invariant $\Psi_X^{u_1}$ for closed manifolds $X$, we need Dirac operators on $\hat{X}_1$ to define the relative invariant.
Take a spin-c structure $\hat{\fs}_1$ of $\hat{X}_1$ with $c_1(\det \hat{\fs}_1) = -c_1(\hat{Q}_1)$. Let $\hat{E}_1$ be the rank 2 complex vector bundle over $\hat{X}_1$ associated with $\hat{Q}_1$. For any connection $A$ on $\hat{Q}_1$, we have the twisted Dirac operator 
\[
\Dir_A:L_1^2( \hat{S}^+_1 \otimes \hat{E}_1) 
\longrightarrow
L^2(\hat{S}_1^- \otimes \hat{E}_1).
\]
Here $\hat{S}^{\pm}_1$ are the spinor bundle associated with $\hat{\fs}_1$. For a connection $A$ with limit $\rho$, put
\[
i_{\rho} := \ind \Dir_A \in \Z.
\]

For flat connections $\rho$ over $Y$ with $\delta_{Y}([\rho]) \equiv 2d \mod 8$, we have 0-dimensional moduli spaces $M_{\hat{P}_1, \rho}$. We define
\[
<\psi_{X_1}^{u_1}, [\rho]> :=
\left\{
\begin{array}{cll}
\# M_{ \hat{P}_1, \rho} & \mod 2 & \text{ if $i_{\rho} \equiv 1 \mod 2 $}, \\
0 & \mod 2 & \text{ otherwise }.
\end{array}
\right.
\]
These numbers define an element $\psi_{X_1}^{u_1} \in C_{d}(Y)$.

\begin{lem}
$\partial \psi_{X_1}^{u_1} = 0$.
\end{lem}

This is proved by counting the number of the ends of 1-dimensional moduli spaces. The proof is similar to that of Lemma \ref{lem bd bd} and we omit the proof. 

\begin{dfn}
$\Psi_{X_1}^{u_1} = [\psi_{X_1}^{u_1}] \in I_{d}(Y)$.
\end{dfn}

Next consider a class $[\Sigma_1] \in H_2(X_1;\Z)$ represented by a closed, oriented surface $\Sigma_1$ in $X_1$. Suppose that the self-intersection number $[\Sigma_1] \cdot [\Sigma_1]$ is even. For flat connections $\rho$ with $\delta_{Y}([\rho]) \equiv 2d-2 \mod 8$, we have 2-dimensional moduli spaces $M_{ \hat{P}_1, \rho}$. By the index theorem and the assumption $[\Sigma_1] \cdot [\Sigma_1] \equiv 0 \mod 2$ we can see that the numerical index of the twisted $\bar{\partial}$ operators over $\Sigma_1$ is even. This implies that the determinant line bundle $\tilde{\cL}_{\hat{\Sigma_1}}$ over the framed moduli space $\tilde{M}_{ \hat{P}_1, \rho}$ descends to the line bundle $\cL_{\hat{\Sigma}_1}$ over $M_{ \hat{P}_1, \rho}$ as in subsection \ref{ss FF}. As in \cite{S torsion}, we can take a section $s_{\Sigma}$ such that the zero locus $M_{\hat{P}_1, \rho} \cap s_{\Sigma_1}^{-1}(0)$ is compact, smooth manifold of dimension 0, i.e., a finite set. We put
\[
< \psi_{X_1}^{u_1}([\Sigma_1]), [\rho] > :=
\left\{
\begin{array}{cll}
\# M_{X_1, \rho} \cap s_{\Sigma_1}^{-1}(0) & \mod 2 & \text{ if $i_{\rho} \equiv 1 \mod 2$, } \\
0 & \mod 2 & \text{ otherwise. }
\end{array}
\right.
\]
These numbers define the element $\psi_{X_1}^{u_1}([\Sigma_1]) \in C_{d-1}(Y)$. As before, this element is a cycle and gives an element of $I_{d-1}(Y)$.

\begin{dfn}
$\Psi_{X_1}^{u_1}([\Sigma_1]) = [ \psi_{X_1}([\Sigma_1])] \in I_{d-1}(Y)$.
\end{dfn}

Lastly consider a relative homology class $[\Sigma_1] \in H_2(X_1,Y;\Z)$. Here $\Sigma_1$ is a compact oriented surface in $X_1$ with boundary $\gamma$, and $\gamma$ is a simple closed curve as in the previous section. Suppose that
\begin{eqnarray}
&& < c_1(Q_1; \varphi), [\Sigma_1] > \equiv 1 \mod 2 \ \text{and} \
[\gamma] \not= 0 \text{ in $H_1(Y;\Z)$, or} 
\label{eq c1 1}
\\
&& < c_1(Q_1;\varphi), [\Sigma_1] > \equiv 0 \mod 2 \ \text{and} \
[\gamma] = 0 \text{ in $H_1(Y;\Z)$}.
\label{eq c1 2}
\end{eqnarray}
Here $c_1(Q_1;\varphi) \in H^2(X_1,Y;\Z)$ is the relative Chern class of $Q_1$  defined by the fixed trivialization $\varphi$ over $Y$.
The conditions above will be needed to obtain the determinant line bundle $\cL_{\hat{\Sigma}_1}$ over $M_{ \hat{X}_1, \rho}$.  It is easy to see that $< c_1(Q_1;\varphi);[\Sigma_1] >$ is independent of the choices of $\varphi$.

\vspace{2mm}

First we consider the case when $d \not \equiv 0 \mod 4$. For generators $[\rho ] \in CF_{2d}(Y) \subset C_{d}(Y;\gamma)$, we have 0-dimensional moduli spaces $M_{\hat{P}_1, \rho}$. Put
\[
< \psi_{X_1}^{u_1}([\Sigma_1]), [\rho] > :=
\left\{
\begin{array}{cll}
\# M_{ \hat{P}_1, \rho } & \mod 2 & \text{ if $i_{\rho} \equiv 1 \mod 2$ }, \\
0 & \mod 2 & \text{ otherwise }.
\end{array}
\right.
\]
For generators $[\rho] \in CF_{2d-2}(Y) \subset C_{d}(Y;\gamma)$, we have 2-dimensional moduli spaces $M_{ \hat{P}_1, \rho }$. We want to define $< \psi_{X_1}([\Sigma_1]), [\rho]> \in \Z_2$ using the determinant line bundle of $\hat{\Sigma}_1$ over $M_{\hat{P}_1, \rho}$. To do this, we must check that the determinant line bundle over $\tilde{M}_{ \hat{P}_1, \rho}$ descends to a line bundle over $M_{ \hat{P}_1, \rho}$ as usual. Here $\tilde{M}_{\hat{P}_1, \rho}$ is the quotient of the space of instantons over $\hat{X}$ by the group of gauge transformations with limit $1$. It is sufficient to show that the index $\ind^{+} \bar{\partial}_{A}$ of twisted $\bar{\partial}$-operator over $\hat{\Sigma}_1$ is even for connections $A$ with limit $\rho$.

\begin{lem} \label{lem ind c1}
Let $\eta$ be the trivial flat connection on $Q_1|_Y$ with respect to the fixed trivialization $\varphi$. For connections $\tilde{A}$ on $\hat{Q}_1$ with limit $\eta$, we have
$\ind^{+} \bar{\partial}_{ \tilde{A} } \equiv < c_1(Q_1;\varphi), [\Sigma_1]> \mod 2$.
\end{lem}

The proof of this lemma will be given at the end of this subsection. By this lemma and Lemma \ref{lem ind even}, we obtain

\begin{lem}
Let $\rho$ be a flat connection on $P_1|_Y$ with $\Gamma_{\rho} \cong U(1)$ and $A$ be a connection on $\hat{P}_1$ with limit $\rho$. We denote by $\tilde{A}$ the lift of $A$ to $\hat{Q}_1$ with the fixed determinant $a_{\det}$.
Under the condition (\ref{eq c1 1}) or (\ref{eq c1 2}), the index $\ind^+ \bar{\partial}_{ \tilde{A} }$ is even.
\end{lem}

By this lemma, we have the determinant line bundle $\cL_{ \hat{\Sigma}_1 } \rightarrow M_{ \hat{P}_1, \rho}$, provided that (\ref{eq c1 1}) or (\ref{eq c1 2}) holds.

For generators $[\rho] \in CF_{2d}(Y) \subset C_{d}(Y;\gamma)$, we can take a section $s_{\hat{\Sigma}_1}$ of $\cL_{\hat{\Sigma}_1} \rightarrow M_{\hat{P}_1, \rho}$ compatible with gluing maps as before. The zero locus $M_{ \hat{P}_1, \rho} \cap s_{\hat{\Sigma}_1}^{-1}(0)$ is a compact smooth 0-dimensional manifold, i.e. a finite set. We define
\[
< \psi_{X_1}^{u_1}([\Sigma_1]), [\rho]> :=
\left\{
\begin{array}{cll}
\# M_{ \hat{P}_1, \rho} \cap s_{ \hat{\Sigma}_1 }^{-1}(0)  & \mod 2 
& \text{ if $i_{\rho} \equiv 1 \mod 2$, } \\
0 & \mod 2 & \text{ otherwise.}
\end{array}
\right.
\]
Put $\psi_{X_1}^{u_1}([\Sigma_1]) := \sum < \psi_{X_1}^{u_1}([\Sigma_1]),[\rho]> [\rho] \in C_{d}(Y;\gamma)$. As usual we have:

\begin{lem}
$\bdg(\psi_{X_1}^{u_1}([\Sigma_1])) = 0$.
\end{lem}

\begin{dfn}
$\Psi_{X_1}^{u_1}([\Sigma_1]) := [\psi_{X_1}^{u_1}([\Sigma_1])] \in I_{d}(Y;\gamma)$.
\end{dfn}

Next we consider the case $d \equiv 0 \mod 4$. Continuously we suppose (\ref{eq c1 1}) or (\ref{eq c1 2}) holds. The only difference from the previous case is the term of the trivial connections. ¡¡We have the 0-dimensional moduli space $M_{\hat{P}_1, \theta}$. We define
\begin{equation} \label{eq psi theta}
\begin{split}
&< \psi_{X_1}^{u_1}([\Sigma_1]), [\theta]> := \\
&\left\{
\begin{array}{cll}
\# M_{\hat{P}_1, \theta} & \mod 2 & \text{ if $i_{\theta} \equiv 1 \mod 2$, $<c_1(Q_1; \varphi),[\Sigma_1]> \equiv 1 \mod 2$, } \\
0 & \mod 2 & \text{ otherwise. } 
\end{array}
\right.
\end{split}
\end{equation}
The other terms are defined as before. We can show that $\psi^{u_1}_{X_1}([\Sigma_1]) \in C_0(Y;\gamma)$ is a cycle and we obtain the relative invariant $\Psi^{u_1}_{X_1}([\Sigma_1]) \in I_0(Y;\gamma)$.

\vspace{2mm}

\noindent
{\it Proof of Lemma \ref{lem ind c1}}

To prove Lemma \ref{lem ind c1}, take a compact, oriented surface $\Sigma_2$ with boundary $\gamma$. Using the restriction $\varphi|_{\gamma}$ of the trivialization, extend $Q_1|_{ \Sigma_1}$ to $\Sigma = \Sigma_1 \cup_{\gamma} \Sigma_2$. We denote it by $Q_{\Sigma}$. Let $\eta_{2}$ be the trivial connection on the trivial $U(2)$-bundle over $\hat{\Sigma}_2 = \Sigma_2 \cup ( \gamma \times \R_{ \geq 0})$. Since $\bar{\partial}_{\eta_2}$ is the direct sum of  two copies of the usual $\bar{\partial}$-operator, $\ind^- \bar{\partial}_{\eta_2}$ is even. Hence
\[
\ind^+ \dbar_{A_{\eta}}  \equiv 
\ind^+ \dbar_{A_{\eta}} + \ind^- \dbar_{ \eta_2} \mod 2.
\]
Here $A_{\eta}$ is a connection on $\hat{Q}_1$ with limit $\eta$.
Moreover we have
\[
\ind^+ \dbar_{A_{\eta}} + \ind^- \dbar_{\eta_2} = 
\ind \dbar_{ A_{\Sigma} }.
\]
Here $A_{\Sigma}$ is the connection over $\Sigma$ obtained by gluing $A_{\eta}$ and $\theta_2$. By index formula we have
\[
\ind \dbar_{ A_{\Sigma} } 
= <c_1(Q_{\Sigma}), [\Sigma] >.
\]
The right hand side is equal to $< c_1(Q_1;\varphi), [\Sigma_1]>$, since there is no contribution from $\Sigma_2$. Thus we have obtained 
\[
\ind^{+} \dbar_{A_{\theta}} \equiv < c_1(Q_1; \varphi), [\Sigma_1]> \mod 2
\]
as required.

\subsection{Gluing formula}

In this subsection, we construct gluing formulas for $\Psi_{X}^{u_1}$. To do this, we need  pairings on Floer homologies. 

\begin{lem}
Let $\rho$ be a flat connection over $Y$ with $\Gamma_{\rho} \cong U(1)$. Then we have
\[
\delta_{-Y}([\rho]) \equiv - \delta_Y([\rho]) - 2 \mod 8.
\]
\end{lem}
The proof is standard and we omit the proof. By this lemma, we have the natural pairing
\[
< \cdot, \cdot >: CF_{2i}(Y) \otimes CF_{-2i-2}(-Y) \rightarrow \Z_2.
\]
This paring induces the pairings
\[
C_{i}(Y) \otimes C_{-i-1}(-Y) \rightarrow \Z_2, \quad
C_{i}(Y; \gamma) \otimes C_{-i}(-Y; \gamma) \rightarrow \Z_2,
\]
which give the identifications $C_{i}(Y)^* = C_{-i-1}(-Y)$, $C_{i}(Y;\gamma)^* = CFF_{-i}(-Y;\gamma)$.
\[
\xymatrix@R=6pt{
C_{i}(Y,\gamma) \ar[r]^{ \hspace{-3mm} \text{dual} } \ar@{=}[d] & C_{-i}(-Y,\gamma) \ar[l] \ar@{=}[d] & 
C_0(Y,\gamma) \ar@{=}[d] \ar[r]^{ \hspace{-2mm} \text{dual} } & C_0(-Y,\gamma) \ar@{=}[d] \ar[l] \\
CF_{2i}(Y) \ar[rdd] & CF_{-2i}(-Y) \ar[ldd] & 
CF_0(Y) \ar[rdd] & CF_0(-Y) \ar[ldd] \\
\oplus & \oplus & \oplus & \oplus \\
CF_{2i-2}(Y) \ar[ruu] & CF_{-2i-2}(-Y) \ar[luu] & 
CF_{-2}(Y) \ar[ruu] & CF_{-2}(-Y) \ar[luu] \\
& & \oplus & \oplus \\
& & \Z_2<[\theta]> \ar[r] & \Z_2<[\theta]> \ar[l]
}
\]
It is easy to see that the pairings induce pairings 
\[
< \ , \ >:I_{i}(Y) \otimes I_{-i-1}(-Y) \rightarrow \Z_2, 
\quad
< \ , \ >:I_{i}(Y;\gamma) \otimes I_{-i}(-Y;\gamma) \rightarrow \Z_2.
\]
Let $X$ be a simply connected, non-spin, closed 4-manifold with a decomposition $X = X_1 \cup_Y X_2$. Here $X_1$ and $X_2$ are simply connected, non-spin 4-manifolds with $b^+ > 0$ and with boundaries $Y$ and $-Y$ respectively. Take a homology class $[\Sigma] \in H_2(X;\Z)$ with $[\Sigma] \cdot [\Sigma] \equiv 0 \mod 2$. Here $\Sigma$ is an embedded surface in $X$. 

\begin{thm} \label{thm gluing 1}
If $\Sigma \subset X_1$ or $\Sigma \subset X_2$, then
\[
\Psi_{X}^{u_1}([\Sigma]) = < \Psi_{X_1}^{u_1}([\Sigma]), \Psi_{X_2}^{u_1} >
\text{or} \ 
\Psi_{X}^{u_1}([\Sigma]) = < \Psi_{X_1}^{u_1}, \Psi_{X_2}^{u_1}([\Sigma]) >.
\]
\end{thm}

Suppose that $\Sigma$ and $Y$ intersect transversely and the intersection $\gamma := Y \cap \Sigma$ is diffeomorphic to $S^1$. We denote $\Sigma \cap X_1$ and $\Sigma \cap X_2$ by $\Sigma_1$ and $\Sigma_2$ respectively. Assume that $[\Sigma_1]$ satisfies (\ref{eq c1 1}) or (\ref{eq c1 2}). ( We can easily see that $[\Sigma_2]$ also satisfies (\ref{eq c1 1}) or (\ref{eq c1 2}).) We have the relative invariants $\Psi_{X_1}^{u_1}([\Sigma_1]), \Psi_{X_2}^{u_1}([\Sigma_2])$.

\begin{thm} \label{thm gluing 2}
Under the above situation, 
\[
\Psi_{X}^{u_1}([\Sigma]) = 
< \Psi_{X_1}^{u_1}([\Sigma_1]), \Psi_{X_2}^{u_1}([\Sigma_2]) >.
\]
\end{thm}

We give outline of the proof of Theorem \ref{thm gluing 2} in the case when $d \equiv 0 \mod 4$. 

Suppose that $d \equiv 0 \mod 4$. Take a sequence $\{ T^{\alpha} \}_{\alpha}$ of positive numbers with $T^{\alpha} \rightarrow \infty$ and a sequence $\{ g^{\alpha} \}_{\alpha}$ of Riemannian metrics on $X$ such that a neighborhood of $Y$ in $X$ is isometric to $(Y \times [-T^{\alpha}, T^{ \alpha }], g_Y + dt^2)$. Here $g_Y$ is the standard metric on $Y$ and $t$ is the coordinate of $[-T^{\alpha}, T^{\alpha}]$. Let $M_P(g^{\alpha})$ be the moduli space of instantons over the Riemannian manifold $(X, g^{\alpha})$ of dimension $3$, where $P$ is an $SO(3)$-bundle with $w_2(P) = w_2(X)$. 
Take sections $s_{\Sigma}^{\alpha}:M_{P}(g^{\alpha}) \rightarrow \cL_{\Sigma}$ compatible with the gluing maps as usual. Then we have

\begin{lem}
Any sequence $[A^{\alpha}] \in M_P(g^{\alpha}) \cap (s_{\Sigma}^{\alpha})^{-1}(0)$ has a subsequence $[A^{\alpha'}]$ such that
\[
[A^{\alpha'}] \longrightarrow ([A_1], [A_2]),
\]
and one of the following occurs:

\begin{enumerate}
\item
$[A_1] \in M_{\hat{P}_1, \rho} \cap s_{ \hat{\Sigma}_1 }^{-1}(0)$, 
$[A_2] \in M_{ \hat{P}_2, \rho}$,
$\Gamma_{\rho} \cong U(1)$, 
$\dim M_{ \hat{P}_1, \rho} = 2, \dim M_{\hat{P}_2, \rho} = 0$ 
(i.e. $[\rho] \in CF_{-2}(Y)$).

\item
$[A_1] \in M_{\hat{P_1}, \rho}$, 
$[A_2] \in M_{ \hat{P}_2, \rho} \cap s_{\hat{\Sigma}_2}^{-1}(0)$,
$\Gamma_{\rho} \cong U(1)$, 
$\dim M_{ \hat{P}_1, \rho} = 0, \dim M_{\hat{P}_2, \rho} = 2$
(i.e. $[\rho] \in CF_0(Y)$).

\item
$[A_1] \in M_{ \hat{P}_1, \theta}$,
$[A_2] \in M_{ \hat{P}_2, \theta}$,
$\dim M_{ \hat{P}_1, \theta} = \dim M_{ \hat{P}_2, \theta} = 0$.

\end{enumerate}

\end{lem}

Take generic sections $s_{\Lambda}^{\alpha}:M_{P}(g^{\alpha}) \rightarrow \Lambda$.
For $\alpha$ sufficiently large, $M_{ P }(g^{\alpha}) \cap (s_{\Sigma}^{ \alpha })^{-1}(0) \cap (s_{\Lambda}^{ \alpha })^{-1}(0)$ is identified with
\[
\begin{split}
&  \coprod_{ 
\begin{subarray}{c}
[\rho]: \\
\delta_{Y}([\rho]) \equiv -2 \mod 8 
\end{subarray}
}
\coprod_{ \fa }
\left( U(1)_{\fa} \cap (s_{\Lambda}^{\alpha})^{-1}(0) \right)
\ \cup \\
& \coprod_{ 
\begin{subarray}{c}
[\rho]: \\
\delta_{Y}([\rho]) \equiv 0 \mod 8 
\end{subarray}
}
\coprod_{ \fb } 
\big( U(1)_{\fb} \cap (s_{\Lambda}^{\alpha})^{-1}(0) \big) 
\ \cup \\
& \coprod_{ \fc } 
SO(3)_{\fc} \cap ( s_{\Sigma}^{\alpha} )^{-1}(0) \cap (s_{\Lambda}^{\alpha})^{-1}(0).
\end{split}
\]
Here $\fa$, $\fb$ and $\fc$ run over
\[
( M_{ \hat{P}_1, \rho} \cap s_{ \hat{\Sigma}_1}^{-1}(0) ) \times M_{ \hat{P}_2, \rho},
\ 
M_{ \hat{P}_1, \rho} \times ( M_{ \hat{P}_2, \rho} \cap s_{ \hat{\Sigma}_2 }^{-1}(0) ),
\ 
\text{and} \ 
M_{ \hat{P}_1, \theta} \times M_{ \hat{P}_2, \theta}
\]
respectively, and $U(1)_{\fa}$, $U(1)_{\fb}$ and $SO(3)_{\fc}$ are the gluing parameters as before. 

\begin{lem} \label{lem gl para l b}
The restrictions $\Lambda|_{ U(1)_{\fa} }$ and $\Lambda|_{ U(1)_{\fb} }$ are non-trivial if and only if $i_{\rho} \equiv 1 \mod 2$,  and $\Lambda|_{ SO(3)_{\fc} }$ and $\cL_{ \Sigma }|_{ SO(3)_{\fc} }$ are non-trivial if and only if $i_{\theta} \equiv 1 \mod 2$ and $< c_1(Q_1;\varphi), [\Sigma_1] > \equiv 1 \mod 2$ respectively.
\end{lem}

This can be proved as \cite[Lemma 3.30]{S Floer}. Note that the condition that $< c_1(Q_1;\varphi), [\Sigma_1] > \equiv 1 \mod 2$ implies that the index of twisted $\dbar$-operator of $\hat{\Sigma}_1$ is odd by Lemma \ref{lem ind c1}. By Lemma \ref{lem gl para l b} we obtain
\begin{gather*}
\Psi_{X}^{u_1}([\Sigma]) 
\equiv  \sum_{ \begin{subarray}{c}
[\rho]: \\
\delta_Y([\rho]) \equiv -2 \mod 8 \\
i_{\rho} \equiv 1 \mod 2
\end{subarray}
} 
\# \left( M_{ \hat{P}_1, \rho} \cap s_{ \hat{\Sigma}_1 }^{-1}(0) \right) \cdot
\# M_{ \hat{P}_2, \rho}
+ \\
\sum_{ \begin{subarray}{c}
[\rho]: \\
\delta_Y([\rho]) \equiv 0 \mod 8 \\
i_{\rho} \equiv 1 \mod 2
\end{subarray}
}
\# M_{ \hat{P}_1, \rho} \cdot 
\# \left( M_{ \hat{P}_2, \rho} \cap s_{ \hat{\Sigma}_2}^{-1}(0)  \right) + \\
\left\{
\begin{array}{cll}
\# M_{ \hat{P}_1, \theta} \cdot \# M_{ \hat{P}_2, \theta} & \mod 2 
& \text{ if $i_{\theta} \equiv 1, < c_1(Q_1; \varphi), [\Sigma_1]> \equiv 1 \mod 2$ }, \\
0 & \mod 2 & \text{ otherwise }.
\end{array}
\right.
\end{gather*}
The right hand side is $< \Psi_{X_1}^{u_1}([\Sigma_1]), \Psi_{ X_2 }^{u_1}([\Sigma_2])>$ by definition. Thus we have proved Theorem \ref{thm gluing 2}.

\section{Calculation and Application}

\subsection{Moduli space} \label{ss moduli space}
In this section, we calculate Floer homology, making use of the results of Austin \cite{Au} and Furuta-Hashimoto \cite{FH}, \cite{Fur inv inst}. 

Throughout this section we assume that $p$ is an odd positive integer. 
Instantons on $L(p,q) \times \R$ correspond to $\Z_p$-invariant instantons on $S^4$. Let $\tilde{P}$ be the  principal $SU(2)$-bundle over $S^4$ with $c_2 = k$ and $M_k$ be the moduli space of instantons on $\tilde{P}$. The moduli space of instantons over $L(p,q) \times \R$ is identified with the fixed point set of a $\Z_p$-action on $M_k$. 

First we consider the action of $T =S^1 \times S^1$ on $S^4 = \C^2 \cup \{ \infty \}$ defined by
\[
(t_1, t_2) \cdot (z_1, z_2) = (t_1 z_1, t_2 z_2).
\]
The set of the isomorphism classes of $SO(3)$-bundles $P$ over $S^4$ with a lift of the $T^2$-action and with $p_1(P) < 0$ is isomorphic to $\Z_{ > 0} \times \Z_{ > 0 }$. We denote the bundle corresponding to $(k_1, k_2) \in \Z_{ > 0} \times \Z_{ > 0}$ by $P(k_1, k_2)$. The bundle is characterized by the following:

\begin{enumerate}[(i)]
\item
The isotropy representation of $t=(t_1, t_2) \in T$ at $\infty \in S^4$ is
\[
\left(
\begin{array}{ccc}
\cos \theta & -\sin \theta & 0 \\
\sin \theta & \cos \theta & 0 \\
0 & 0 & 1
\end{array}
\right)
\]
up to conjugate. Here $t_1^{k_1} t_2^{k_2} = e^{ i \theta }$.

\vspace{2mm}

\item
The isotropy representation of $t \in T$ at $0 \in S^4$ is
\[
\left(
\begin{array}{ccc}
\cos \theta & -\sin \theta & 0 \\
\sin \theta & \cos \theta & 0 \\
0 & 0 & 1
\end{array}
\right)\]
up to conjugate. Here $t_1^{k_1} t_2^{-k_2} = e^{i \theta}$.

\vspace{2mm}

\item
$p_1(P(k_1, k_2)) = -4k_1 k_2$.

\end{enumerate}

\vspace{2mm}

Let $\tilde{P}(k_1, k_2)$ be the $SU(2)$-bundle with $\tilde{P}(k_1, k_2)/ \{ \pm 1 \} = P(k_1,k_2)$. The second Chern class of $\tilde{P}(k_1,k_2)$ is $k_1 k_2$. A double cover $\tilde{T}$ of $T$ naturally acts on $\tilde{P}$ and we have the induced action of $\tilde{T}$ on the moduli space $M_k = M( \tilde{P}(k_1,k_2))$. Here $k=k_1 k_2$.

By Atiyah-Bott-Lefschetz fixed point formula, we obtain:

\begin{lem} [\cite{Fur inv inst}]
Let $[A] \in M_k$ be a $\tilde{T}$-invariant instanton and $\tilde{t} \in \tilde{T}$. Then we have
\[
\Tr (  \tilde{t}  |  T_{[A]} M ( \tilde{P}(k_1, k_2))  ) =
-1 + \sum_{i, j} a_{ij} t_1^{i} t_2^{j}.
\]
Here $t=(t_1, t_2) \in T$ is the image of $\tilde{t} \in \tilde{T}$ under the projection and
\[
a_{ij} =
\left\{
\begin{array}{cl}
2 & \text{if $|i| < k_1, |j| < k_2$}, \\
1 & \text{if $|i| = k_1, |j| < k_2$ or $|i| < k_1, |j|=k_2$,  } \\
0 & \text{otherwise}.
\end{array}
\right.
\]
\end{lem}

We have a natural inclusion $\Z_p \hookrightarrow T$ defined by $\zeta \mapsto (\zeta, \zeta^q)$. Since we assumed $p$ is odd, there is a unique lift $\Z_p \hookrightarrow \tilde{T}$. Restricting the above formula to $\Z_p \subset \tilde{T}$, we obtain the following:

\begin{cor} \label{cor moduli dim}
The dimension of the fixed points set $M( \tilde{P}(k_1, k_2) )^{\Z_p}$ is given by
\[
\dim M( \tilde{P}(k_1,k_2) )^{ \Z_p } = -1 + 2N_1(k_1,k_2;p,q) + N_2(k_1,k_2;p,q).
\]
Here
\begin{equation} \label{eq N}
\begin{split}
&N_1(k_1,k_2;p,q) = \\
&\quad \# \{ \ (i, j) \in \Z^2 \ | \  i + qj \equiv 0 \mod p,\ |i| < k_1, |j| < k_2  \ \}, \\
&N_2(k_1,k_2;p,q) = \\
&\quad \# \{ \ (i,j) \in \Z^2 \ |  
\text{ $i + qj \equiv 0 \mod p$, $|i| = k_1, |j| < k_2$, or $|i| < k_1, |j|=k_2$  } \}.
\end{split}
\end{equation}
\end{cor}

For flat connections $\rho$ over $Y = L(p,q)$ with $\Gamma_{\rho} \cong U(1)$, 
\[
\delta_Y([\rho]) \equiv \dim M_{\rho \theta} + 1 \mod 8.
\]
On the other hand, the dimension of the moduli space $M_{\rho \theta}$ is congruent to $\dim M( \tilde{P}(k_1,k_2))^{ \Z_p }$ modulo 8. Here $k_1, k_2$ are positive integers such that the restriction of the isotropy representation of $\tilde{P}(k_1, k_2)$ at $\infty \in S^4$ to $\Z_p$ is isomorphic to the holonomy representation of $\rho$ and that at $0 \in S^4$ is trivial. We can find such $k_1, k_2$ as follows. Suppose that the holonomy representation of $\rho$ is given by
\[
1 \longmapsto 
\left(
\begin{array}{cc}
\zeta^l & 0 \\
0 & \zeta^{-l}
\end{array}
\right)
\]
up to conjugation. Here $l$ is a positive integer with $0 < l < p$. Take a positive integer $r$ with $q r \equiv 1 \mod p$. Then $k_1, k_2$ are any positive integers satisfying
\[
k_1 \equiv l \mod p, \quad k_2 \equiv -rl \mod p.
\]

By Corollary \ref{cor moduli dim}, we have:

\begin{cor}
Take  a flat connection $\rho$ over $Y = L(p,q)$ with $\Gamma_{\rho} \cong U(1)$.
Then we have
\[
\delta_Y([\rho]) \equiv 2N_1(k_1, k_2;p,q) + N_2(k_1,k_2;p,q) \mod 8,
\]
where $k_1, k_2 > 0$ are determined as above.
\end{cor}

If $(i, j)$ is a solution to the equation defining $N_2(k_1, k_2;p,q)$, then $(-i, -j)$ is also a solution and $(i, j) \not= (0, 0)$. Therefore $N_2(k_1,k_2;p,q)$ is even. Thus we have:

\begin{cor} \label{cor even}
When $p$ is odd, $\delta_Y([\rho]) \equiv 0 \mod 2$.
\end{cor}

The boundary map of Floer homology $I_*(L(p,q))$ was defined using the moduli spaces $M_{\rho \sigma}$ of dimension $1$. Such moduli spaces are completely determined as follows:

\begin{thm} [\cite{Au},\cite{FH}, \cite{Fur inv inst}] \label{thm 1 dim moduli}

\begin{enumerate}
\item
For any $k_1, k_2 > 0$, $M( \tilde{P}(k_1, k_2))^{ \tilde{T} } = \R_{>0}$.

\vspace{2mm}

\item
Let $\tilde{P} \rightarrow S^4$ be an $SU(2)$-bundle with $c_2 = k$. Suppose that the action of $\Z_p$ on $S^4$ lifts to an action on $\tilde{P}$. 
If the fixed point set $M_k^{ \Z_p }$ is not empty and $1$-dimensional, then there exists $k_1, k_2 > 0$ such that the action of $\Z_p$ on $\tilde{P}$ is the restriction of the action of $\tilde{T}$ on $\tilde{P}(k_1, k_2)$. Furthermore we have an identification $M_k^{\Z_p} = M( \tilde{P}(k_1, k_2))^{ \tilde{T}}$.

\end{enumerate}
\end{thm}

\begin{cor}
Let $\rho, \sigma$ be flat connections over $Y = L(p, q)$ such that the formal dimension of $M_{\rho \sigma}$ is $1$. If there exists $k_1, k_2 > 0$ such that
\begin{enumerate}[(i)]
\item
the isotropy representation of $\tilde{P}(k_1, k_2)$ at $0 \in S^4$ is isomorphic to the holonomy representation of $\rho$,

\item
the isotropy representation of $\tilde{P}(k_1, k_2)$ at $\infty \in S^4$ is isomorphic to the holonomy representation of $\sigma$, and

\item
$\dim M( \tilde{P}(k_1, k_2))^{\Z_p} = 1$,
\end{enumerate}

then we have an identification $M_{\rho \sigma} = \R$. If such $k_1$ and $k_2$ do not exist, then $M_{\rho \sigma} = \emptyset$.
\end{cor}

\subsection{Index of Dirac operator}

Take flat connections $\rho$ and  $\sigma$ on the trivial $SU(2)$-bundle $Q = Y \times SU(2)$ with $\Gamma_{\rho}, \Gamma_{\sigma} \cong U(1)$. Assume that $M_{\rho \sigma}$ is not empty and that $M_{\rho \sigma}$ is 1-dimensional. We will compute the index $\ind \Dir_A$ of the twisted Dirac operator. We write $A'$ for the pull-back of $A$ by the projection $S^3 \times \R \rightarrow Y \times \R$. We have the virtual representation space $\Ind \Dir_{A'}$ of $\Z_p$. We can write
\[
\Ind \Dir_{A'} = \sum_n b_n \chi_n,
\]
where $\chi_n$ is the 1-dimensional representation space of $\Z_p$ of weight $n$. Then we have
\[
\ind \Dir_A = b_0.
\]
For the lifts $\rho', \sigma'$ of $\rho, \sigma$, we may take trivializations $\varphi_1, \varphi_2$ of the trivial $SU(2)$-bundle $Q' = S^3 \times SU(2)$  such that $\rho', \sigma'$ are trivial with respect to $\varphi_1, \varphi_2$. Using $\varphi_1$ and  $\varphi_2$, we extend $\pi^* Q'$ to $S^4 = D^4 \cup S^3 \times \R \cup D^4$. Here $\pi$ is the projection $S^3 \times \R \rightarrow S^3$.  We have the extension $B$ of $A'$ to $S^4$. That is, $B$ is equal to $A'$ over $S^3 \times \R$ and trivial on the Discs.

\begin{lem} \label{lem add rep}
We have
\[
\Ind \Dir_B = \Ind \Dir_{ \theta_{D^4}} + \Ind \Dir_{A'} + \Ind \Dir_{ \theta_{D^4}}
\]
as virtual representation spaces of $\Z_p$. 
\end{lem}

This lemma will be proved later.
\vspace{2mm} 

By the Weizenb\"ock formula and the facts that $D^4$ has a metric of positive scalar curvature which restrict to the standard metric on $\partial D^4 = S^3$ and $\theta_{D^4}$ is flat, we have
\[
\ind \Dir_{\theta_{D^4}} = 0.
\]
Therefore we obtain
\[
\Ind \Dir_{A'} = \Ind \Dir_{B}
\]
as representation spaces of $\Z_p$.
Theorem \ref{thm 1 dim moduli} implies that $B$ is a $\tilde{T}$-invariant connection on $\tilde{P}(k_1, k_2)$ for some $(k_1, k_2) \in \Z_{>0} \times \Z_{>0}$ and we can regard $\Ind \Dir_{B}$ as a representation space of $\tilde{T}$. By Atiyah-Bott-Lefschetz fixed point formula we obtain
\[
\begin{split}
\ind ( \Dir_B, \tilde{t}) 
&= \frac{ 
t_1^{ \frac{k_1}{2} } t_2^{ -\frac{k_2}{2} } + t_1^{ -\frac{k_1}{2} } t_2^{ \frac{k_2}{2} } - t_1^{ \frac{k_1}{2} } t_2^{ \frac{k_2}{2} } - t_1^{ - \frac{k_1}{2} } t_2^{ -\frac{k_2}{2} } }{ (t_1^{ \frac{1}{2} } - t_1^{ - \frac{1}{2} }) 
(t_2^{ \frac{1}{2} } - t_2^{ - \frac{1}{2} }) 
} \\
&=
- ( t_1^{ \frac{ -k_1 + 1 }{2} }  + t_1^{ \frac{ -k_1 + 3 }{2} } + \cdots + t_1^{ \frac{ k_1 - 1 }{2} })
( t_2^{ \frac{ -k_2 + 1 }{2} }  + t_2^{ \frac{ -k_2 + 3 }{2} } + \cdots + t_2^{ \frac{ k_2 - 1 }{2} }).
\end{split}
\]
Restricting to $\Z_p \subset \tilde{T}$, we have
\[
\begin{split}
&\ind (\Dir_B, \zeta) \\
& \quad = - ( \zeta^{ \frac{ -k_1 + 1 }{2} } + \zeta^{ \frac{-k_1 + 3}{2} } + \cdots + \zeta^{ \frac{ k_1 - 1 }{2} })
( \zeta^{  \frac{q(-k_2 + 1)}{2}  } + \zeta^{ \frac{q(-k_2 + 3)}{2} } + \cdots \zeta^{ \frac{q(k_2 - 1)}{2} }) \\
& \quad = - \sum_{ a = 0 }^{ k_1 - 1 } \sum_{ b = 0 }^{ k_2 - 1 }
\zeta^{ \frac{ -k_1 + 2a + 1 + q( -k_2 + 2b + 1 ) }{2} }
\end{split}
\]
The index $\ind \Dir_{A} \in \Z$ is the constant term of the right hand side. Thus we have obtained:

\begin{prop}
The index $\ind \Dir_{A} \in \Z$ is equal to minus the number of solutions of the following equation for $(a, b)$:
\[
-k_1 + 2a + 1 + q(-k_2 + 2b + 1) \equiv 0 \mod 2p
\quad ( 0 \leq a \leq k_1 - 1, 0 \leq b \leq k_2 - 1).
\]
\end{prop}

It remains to determine $(k_1, k_2)$. Suppose that the holonomy representations of $\rho, \sigma$ are given by
\[
\left(
\begin{array}{cc}
\zeta^l & 0 \\
0 & \zeta^{-l}
\end{array}
\right), 
\quad
\left(
\begin{array}{cc}
\zeta^m & 0 \\
0 & \zeta^{-m}
\end{array}
\right)
\]
where $0 < l < p, 0 < m < p$.
From the fact that the restriction of the isotropy representations of $\tilde{P}(k_1,k_2)$ to $\Z_p \subset \tilde{T}$ are given by these matrix (up to conjugation), $k_1$ and $k_2$ must satisfy one of the following four equalities:
\[
k_1 \equiv \pm l + \pm m \mod p, \quad k_2 \equiv r( \pm m - \pm l) \mod p.
\]
Here $r$ is a positive integer with $rq \equiv 1 \mod p$. Note that we must consider both of $l$ and $-l$ since the matrixes 
\[
\left(
\begin{array}{cc}
\zeta^l & 0 \\
0 & \zeta^{-l}
\end{array}
\right)
\text{ and }
\left(
\begin{array}{cc}
\zeta^{-l} & 0 \\
0 & \zeta^l
\end{array}
\right)
\]
are conjugate. Similarly for $m$. Since $\dim M_{\rho \sigma}$ is $1$ and $\dim M_{\rho \sigma}$ is given by the formula in Corollary \ref{cor moduli dim}, $k_1$ and $k_2$ also satisfy the condition that the set of the solutions $(i, j)$ to the equation
\begin{equation*}
i + qj \equiv 0 \mod p, \quad |i| \leq k_1, |j| \leq k_2
\end{equation*}
is a subset of $\{(0,0), \pm (k_1, k_2), \pm (k_1, -k_2) \}$. If $M_{\rho \sigma}$ is not empty, we can find such a pair $(k_1,k_2)$ by Theorem \ref{thm 1 dim moduli}. 
\vspace{3mm}

The discussions of the previous subsection and this subsection give us a way to compute $I_*(Y)$. Here we summarize the way to compute $I_*(Y)$. Fix a positive integer $r$ with $qr \equiv 1 \mod p$.

\begin{enumerate}[(i)]
\item $\delta_{Y}([\rho])$.

For an integer $l$ with $0 < l < p$, let $\rho_l$ be a flat connection whose holonomy representation is given by
\[
\left(
\begin{array}{cc}
\zeta^{l} & 0 \\ 
0 & \zeta^{-l}
\end{array}
\right)
\]
where $\zeta = e^{ \frac{2 \pi \sqrt{-1} }{p} }$.  Choose positive integers $k_1, k_2$ with
\begin{equation} \label{eq klr}
k_1 \equiv l \mod p, \quad k_2 \equiv -rl \mod p.
\end{equation}
Consider the equation for $(i, j) \in \Z \times \Z$:
\begin{equation} \label{eq (i,j)}
i + qj \equiv 0 \mod p, \quad |i| \leq k_1, \ |j| \leq k_2.
\end{equation}
Define $N_1(k_1, k_2;p,q)$ to be the number of solutions $(i,j)$ with $|i|< k_1$, $| j | < k_2$, and define $N_2(k_1, k_2;p,q)$ to be the number of solutions $(i, j)$ with $|i| = k_1, |j| < k_2$ or with $|i| < k_1, |j| = k_2$.
Then the degree $\delta_Y([\rho_l])$ is 
\[
2 N_1(k_1,k_2;p,q) + N_2(k_1,k_2;p,q) \mod 8.
\]
The vector space $C_i(Y)$ is spanned by the gauge equivalence classes $[\rho_l]$ with $\delta_{Y}([\rho_l]) \equiv 2i \mod 8$.

\vspace{2mm}

\item $< \partial([\rho]), [\sigma]>$.

Take generators $[\rho_l] \in C_i(Y)$, $[\rho_m] \in C_{i-1}(Y)$. Here $0 < l < p, \ 0 < m < p$. 

\vspace{2mm}

\begin{enumerate}
\item
If there exists $k_1, k_2 > 0$ such that $k_1$ and $k_2$ satisfy one of the following four equations
\begin{equation} \label{eq pm krlm}
k_1 \equiv \pm l + \pm m \mod p, \quad k_2 \equiv r(\pm m - \pm l) \mod p,
\end{equation}
and the set of solutions to (\ref{eq (i,j)}) is a subset of $\{ (0,0), \pm (k_1, k_2), \pm (k_1, -k_2) \}$, then 
\begin{gather*} 
< \partial([\rho_l]), [\rho_m]>
 \equiv \nonumber \\
\# \left\{
(a, b) \in \Z^2
\left|
\begin{array}{c}
0 \leq a \leq k_1 - 1, \  0 \leq b \leq k_2 - 1, \\
-k_1 + 2a + 1 + q(-k_2 + 2b + 1) \equiv 0 \mod 2p
\end{array}
\right.
\right\}
\mod 2. 
\end{gather*}

\item
Otherwise  $<\partial([\rho_l], [\rho_m]> \equiv 0 \mod 2$.
\end{enumerate}
\end{enumerate}

\noindent
{\it Proof of Lemma \ref{lem add rep}}.

In the proof of the usual addition property of the index, three data are used. These are cut off functions, stabilizations $S_0:\R^N \rightarrow \Gamma (F_1)$ of differential operators $D:\Gamma(F_0) \rightarrow \Gamma(F_1)$ (i.e. $D \oplus S_0$ is surjective) and right inverses of the stabilized operators $D \oplus S_0$. (See \cite{D Floer} for details.) In our setting, $F_0, F_1$ and $D$ are defined over$S^3 \times \R$ or $D^4$. Moreover $\Z_p$ acts on $F_0, F_1$ and $D$ is $\Z_p$-equivalent. It is sufficient to prove that we can make these data $\Z_p$-equivalent.

It easy to see that we can choose $\Z_p$-invariant cut off functions. To take a $\Z_p$-equivalent stabilization, put $S_j := \zeta^j S_0$ for $j=1, 2,\dots, p-1$. Here $S_0:\C^N \rightarrow \Gamma(F_1)$ is a fixed stabilization of $D$ and $\zeta = \exp( 2\pi \sqrt{-1}/p)$. Define a $\Z_p$-action on $\C^{pN}=\C^{N} \oplus \cdots \oplus \C^N$ by
\[
\zeta \cdot (v_0, v_1,\dots , v_{p-1}) = (v_{p-1}, v_0, v_1, \dots, v_{p-2}).
\]
Then
\[
S:= S_0 \oplus S_1 \oplus \cdots \oplus S_{p-1}:\C^{pN} \rightarrow \Gamma(F_1)
\]
is a stabilization of $D$ and $\Z_p$-equivalent. Choose any right inverse $Q'$ of $D \oplus S$ and define
\[
Q := \frac{1}{p} \sum_{j = 0}^{p-1} \zeta^j Q' \zeta^{-j}.
\]
Then this operator is also a right inverse of $D \oplus S$ and $\Z_p$-equivalent.

\subsection{$L(8N + 1, 2)$}
In this subsection, we compute Floer homology for $Y = L(8N+1,2)$. The following vanishing of Floer homology was suggested by Yuichi Yamada.

\begin{prop} \label{prop 8N+1}
For all $i$, we have $I_i(L(8N+1,2)) = 0$.
\end{prop}

Put $Y = L(8N+1, 2)$.
The complex $C_*(Y)$ is generated by $[\rho_1],[\rho_2], \cdots, [\rho_{4N}]$. Note that $[\rho_l] = [\rho_{8N+1-l}]$. Put $r=4N+1$. Then $2 \cdot r = 1 \mod 8N+1$.

\vspace{2mm}

\noindent
$\bullet$ $\delta_Y([\rho_l])$.

\noindent
Let $l$ be an odd integer with $1 \leq l \leq 4N$. Put $k_1 = l, k_2 = 4N - \frac{l-1}{2}$. Then (\ref{eq klr}) is satisfied. We consider the equation (\ref{eq (i,j)}). The solutions to the equation
\[
i + 2j = 0 \quad (|i| \leq k_1 = l, \ |j| \leq k_2 = 4N - \frac{l-1}{2})
\]
are 
\[
\begin{split}
(i,j) =
&(0,0), (2, -1), (4,-2), \dots , (l-1, -\frac{l-1}{2}), \\
& -(2, -1), -(4,-2), \dots , -(l-1, -\frac{l-1}{2}).
\end{split}
\]
The solutions to the equations
\begin{equation*}
i+2j = \pm (8N+1) \quad ( |i| \leq l,  |j| \leq 4N - \frac{l-1}{2})
\end{equation*}
are
\[
(i,j) =  \pm (l, 4N - \frac{l-1}{2}).
\]
For any $m \in \Z$ with $|m| \geq 2$, the equations
\[
i + 2j = m(8N+1), \quad
(|i| \leq l, \ |j| \leq 4N - \frac{l-1}{2})
\]
do not have solutions. Hence we have
\[
N_1(k_1, k_2; 8N+1, 2) = l, \quad N_2(k_1, k_2; 8N+1, 2) = 0.
\]
Therefore
\[
\delta_Y([\rho_l]) \equiv 2l \mod 8.
\]

\vspace{2mm}

Next let $l$ be an even integer with $1 \leq l \leq 4N$. Put $k_1 = l$, $k_2 = 8N + 1 - \frac{l}{2}$. Then (\ref{eq klr}) is satisfied. The solutions to the equation
\[
i + 2j = 0 \quad (|i| \leq k_1 = l, \ |j| \leq k_2 = 8N + 1 - \frac{l}{2})
\]
are
\[
\begin{split}
(i,j) = 
&(0, 0), (2, -1), (4,-2), \dots, (l-2, -\frac{l}{2}+1), (l, -\frac{l}{2}), \\
&-(2,-1), -(4,-2), \dots,-(l-2, -\frac{l}{2}+1), -(l, -\frac{l}{2}).
\end{split}
\]
The solutions to the equations
\[
i + 2j = \pm(8N + 1) \quad (|i| \leq l, \ |j| \leq 8N + 1 - \frac{l}{2})
\]
are
\[
\begin{split}
(i,j) = 
&\pm(1, 4N), \pm(3, 4N-1),\dots, \pm(l-1,4N-\frac{l}{2}+1) \\
&\pm(-1, 4N+1), \pm(-3, 4N+2), \cdots, \pm(-l+1, 4N+\frac{l}{2}).
\end{split}
\]
The solutions to the equations
\[
i + 2j = \pm2(8N+1) \quad (|i| \leq l, \ |j| \leq 8N + 1 - \frac{l}{2})
\]
are
\[
(i,j) = \pm(l, 8N  + 1 - \frac{l}{2}).
\]
For $m \in \Z$ with $|m| \geq 3$, the equations
\[
i + 2j = m (8N+1) \quad (|i| \leq l, \ |j| \leq 8N + 1 - \frac{l}{2})
\]
do not have solutions. Hence we have
\[
\begin{split}
N_1(k_1, k_2;8N+2, 2) &= 1 + 2 \cdot \frac{l-2}{2} + 4 \cdot \frac{l}{2} = 3l-1, \\ 
N_2(k_1, k_2;8N+2, 2) &= 2, \\
\delta_{Y}([\rho]) &\equiv 2 (3l-1) + 2 \equiv 6l \mod 8.
\end{split}
\]
Thus we have obtained
\[
\begin{split}
\delta_{Y}([\rho_l]) 
& \equiv
\left\{
\begin{array}{cl}
2l & \text{if $l$ is odd}, \\
6l & \text{if $l$ is even},
\end{array}
\right. \\
C_i(L(8N+1, 2)) 
&= \left\{
\begin{array}{ll}
\Z_2 < [\rho_3], [\rho_7], \dots, [\rho_{4N-1}] > & i \equiv 3 \mod 4, \\
\Z_2 < [\rho_2], [\rho_6], \dots, [\rho_{4N-2}] > & i \equiv 2 \mod 4, \\
\Z_2 < [\rho_1], [\rho_5], \dots, [\rho_{4N-3}] > & i \equiv 1 \mod 4, \\
\Z_2 < [\rho_4], [\rho_8], \dots, [\rho_{4N}]> & i \equiv 0 \mod 4.
\end{array}
\right.
\end{split}
\]

\noindent
$\bullet$ $< \partial([\rho_{4s+4}]), [\rho_{4t+3}]>$ ($0 \leq s, t \leq N-1$).

If $s \not = t$, for $k_1, k_2>0$ satisfying one of the equations (\ref{eq pm krlm}), we have $k_1 \geq 3$. Hence $(2, -1)$ is a solution to the equation (\ref{eq (i,j)}) and $(2, 1)$ is not included in $\{ (0,0), \pm (k_1, k_2), \pm (k_1, \pm k_2) \}$. Therefore we have
\[
< \partial ([\rho_{4s+4}]), [\rho_{4t+3}])> \equiv 0 \mod 2
\quad (s \not= t).
\]
Suppose that $s = t$. Put $k_1 = 1, k_2 = 4N - 4s - 3$. Then we can see that $k_1$ and $k_2$ satisfy one of the equations (\ref{eq pm krlm}):
\[
k_1 \equiv l - m, \quad k_2 \equiv (4N+1)(-m - l) \mod 8N + 1.
\]
The only solution to the equation
\[
i + 2 j \equiv 0 \mod 8N + 1 \quad (|i| \leq 1, |j| \leq 4N - 4s -3)
\]
is $(0, 0)$. Hence $< \partial([\rho_{4s+4}], [\rho_{4s+3}] >$ is the number of solutions to the equation
\[
2 (-4N + 4s + 3 + 2b + 1 ) \equiv 0 \mod 2(8N+1), \quad
(|b| \leq 4N - 4s - 4).
\]
This equation has the unique solution
\[
b = 2N - 2s - 2.
\]
Therefore 
\[
< \partial([\rho_{4s+4}], [\rho_{4s+3}] > \equiv 1 \mod 2.
\]
Hence the boundary map $C_{4i}(Y) \rightarrow C_{4i-1}(Y)$ is isomorphic. It follows from $\partial \circ \partial = 0$ that  $\partial:C_{4i+3}(Y) \rightarrow C_{4i+2}(Y)$ is zero.

\vspace{2mm}

$\bullet$
$< \partial ([\rho_{4s+2}]),[\rho_{4t+1}] >$ $(0 \leq s, t \leq N - 1)$

As in the previous case, we can see that $<\partial ([\rho_{4s+2}]), [\rho_{4t+2}]> \equiv 0 \mod 2$ for $s \not = t$.

Let $s = t$. Put $k_1 = 1, k_2 = 4N - 4s - 1$. Then $k_1$ and $k_2$ satisfies one of the equations (\ref{eq pm krlm}):
\[
k_1 \equiv l - m, \ k_2 \equiv (4N + 1)(-m - l) \mod 8N + 1.
\] 
The equation (\ref{eq (i,j)}) with $k_1 = 1, k_2 = 4N - 4s -1$ has the unique solution $(0,0)$. Moreover the equation 
\[
2(-4N + 4s + 1 + 2b + 1) \equiv 0 \mod 2(8N + 1), \ 
(0 \leq b \leq  4N - 4s - 2)
\]
also has the unique solution
\[
b = 2N - 2s + 1.
\]
Hence we have
\[
< \partial ([\rho_{4s+2}]), [\rho_{4s+1}] > \equiv 1 \mod 2.
\]
The boundary map $\partial:C_{4i + 2}(Y) \rightarrow C_{4i + 1}(Y)$ is isomorphic and $\partial:C_{4i + 1 }(Y) \rightarrow C_{4i}(Y)$ is zero:
\[
\cdots \stackrel{0}{\rightarrow}
C_{4i + 4}(Y) \stackrel{\cong}{\rightarrow}
C_{4i + 3}(Y) \stackrel{0}{\rightarrow}
C_{4i + 2}(Y) \stackrel{\cong}{\rightarrow}
C_{4i + 1}(Y) \stackrel{0}{\rightarrow}
C_{4i}(Y) \stackrel{\cong}{\rightarrow}
\cdots
\]
Thus we have obtained the statement of Proposition \ref{prop 8N+1}.

\subsection{Application}
The aim of this subsection is to prove the following:

\begin{thm} \label{thm dec}
Let $X = \CP^2 \# \CP^2$ and $Y = L(p,q)$. Then $X$ does not admit a decomposition $X = X_1 \cup_{Y} X_2$, when $p$ is a prime number of the form $16N + 1$ and $q = 2$. Here $X_1$ and $X_2$ are simply connected, non-spin 4-manifolds with $b^+ = 1$ and with $\partial X_1 = Y$, $\partial X_2 = -Y$.
\end{thm}

Before proving this, we give some relevant remarks.

\begin{rem}
\begin{enumerate}[(1)]
\item
The first remark was pointed out by Kouichi Yasui. Assume that $X = \CP^2 \# \CP^2$ has a decomposition $X = X_1 \cup_Y X_2$ for some $Y = L(p,q)$. Here $X_1$ and $X_2$ are as in Theorem \ref{thm dec}. Then $p$ must be of the form $a^2 + b^2$ for some integers $a$ and $b$. This can be seen as follows. As stated in Lemma \ref{lem X dec} below, we have $H_2(X_1;\Z) = \Z$. Take a generator $\sigma \in H_2(X_1;\Z)$. Since $H_1(Y;\Z) = \Z /(\sigma^2)$, we have $\sigma^2 = p$. Here $\sigma^2$ is the self-intersection number of $\sigma$. On the other hand, we can think of $\sigma$ as a homology class of $X$ through the natural map $H_2(X_1;\Z) \rightarrow H_2(X;\Z)$. In $H_2(X;\Z)$, we can write
\[
\sigma = a H_1 + b H_2,
\]
where $H_1$ and $H_2$ are the natural generators of $H_2(X;\Z)$, and $a, b$ are some integers. Hence we have $\sigma^2 = a^2 + b^2$.
Therefore we have $p = a^2 + b^2$. 

From this, we see that $X = \CP^2 \# \CP^2$ does not admit a decomposition along $L(p, q)$ for $p = 3, 6, 7, 11, \dots$.

\vspace{2mm}

\item
By Dirichlet's theorem (see, for example, \cite[Chapter VI]{Ser}), there are infinitely many prime numbers of the form $16N + 1$:
\[
17, 97, 113, 193, 241, 257, 337 \dots
\]
By Fermat's two squares theorem, every prime number $p$ with $p \equiv 1 \mod 4$ is a sum of two squares (See \cite{Ded}). Hence all prime numbers of the form $16N + 1$ can be written as $a^2 + b^2$ for some $a, b \in \Z$.

\vspace{2mm}

\item
M. Tange and Y. Yamada \cite{TY} showed that there are decompositions  $\CP^2 \# \CP^2 = X_1 \cup_Y X_2$ for infinitely many lens spaces $Y = L(p,q)$. Here $X_1, X_2$ are as in Theorem \ref{thm dec}. For examples,
\[
\begin{split}
&Y = L(5,1), L(13,9), L(29, 9), L(34,9), \\
& \qquad \qquad L(89,25), L(233, 64), \dots, L(28657, 7921), \dots
\end{split}
\]
On the other hand, $28657$ is a prime number and congruent to $1$ modulo $16$.¡¡(This was pointed out by Yamada.) Hence $\CP^2 \# \CP^2$ can not be decomposed along $L(28657, 2)$ by Theorem \ref{thm dec}.
\end{enumerate}
\end{rem}

It is easy to see the following.

\begin{lem} \label{lem X dec}
Put $Y = L(p,q)$. Assume that $X = \CP^2 \# \CP^2$ has a decomposition $X = X_1 \cup_{Y} X_2$. Here $X_1, X_2$ are as in Theorem \ref{thm dec}. Then we have
\[
H^2(X_i; \Z) = \Z, \quad H^2(X_i;\Z_2) = \Z_2.
\]
\end{lem}

We show the following lemma making use of Proposition \ref{prop moduli X_1 dim}.

\begin{lem} \label{lem X dec 2}
Let $Y$ be $L(p,q)$ and $X_1$ as in Theorem \ref{thm dec}. Take an $SO(3)$-bundle $P_1$ over $X_1$ with $w_2(P_1) = w_2(X_1)$. Then we have
\[
\dim M_{ \hat{P}_1, \rho} \equiv - \delta_{Y}([\rho]) - 2p - 6 \mod 8.
\]
\end{lem}

Take a cohomology class $\alpha \in H^2(X_1;\Z)$ with the properties in Lemma \ref{lem alpha}. From Proposition \ref{prop moduli X_1 dim}, we have only to show that $\alpha^2 \equiv p \mod 8$.
Since $X_1$ is non-spin and $\alpha \equiv w_2(X_1) \mod 2$, we can take $p \beta$ as $\alpha$. Here $\beta \in H^2(X_1;\Z) = \Z$ is a generator. (Recall that we assumed $p$ is odd.) By the exact sequence 
\[
\begin{split}
& H^1(Y;\Z) = 0 \rightarrow H^2(X_1, Y;\Z) \stackrel{j^*}{\rightarrow} \\
& \qquad H^2(X_1;\Z) = \Z \rightarrow H^2(Y;\Z) = \Z_p \rightarrow H^3(X_1,Y;\Z) = 0,
\end{split}
\]
we see that $H^2(X_1,Y;\Z) = \Z$ and that $j^* (\tilde{\alpha}) = p \beta = \alpha$ for some generator $\tilde{\alpha} \in H^2(X_1,Y;\Z) $. Since the pairing
\[
H^2(X_1,Y;\Z) \otimes H^2(X_1;\Z) \longrightarrow \Z
\]
induces an identification $H^2(X_1,Y;\Z) = H^2(X_1;\Z)^*$, we have
\[
< \tilde{\alpha} \cup \beta, [X_1,Y]> = \pm 1.
\]
Thus
\[
\alpha^2 = < \tilde{\alpha} \cup \alpha, [X_1,Y]> = \pm p.
\]
Since $H^2(X_1;\Z) = \Z$ and $b^+(X_1) = 1$, the above equality has a plus sign, and hence $\alpha^2 = p$.

\vspace{3mm}

\noindent
{\it Proof of Theorem \ref{thm dec}}

Put $Y = L(p, 2)$, where $p$ is a prime number of the form $16N+1$. Assume that $X = \CP^2 \# \CP^2$ admits a decomposition $X  = X_1 \cup_Y X_2$. Here $X_1$ and $X_2$ as in Theorem \ref{thm dec}. It follows from Lemma \ref{lem X dec 2} that 
\begin{equation} \label{eq moduli X_1}
\dim M_{X_1, \rho} \equiv -\delta_Y([\rho])  \mod 8.
\end{equation}

By \cite[Theorem 3.29]{S Floer}, we can take a cohomology class $h_0 \in H_2(X;\Z)$ with $h_0 \cdot h_0 \equiv 0 \mod 2$ such that
\[
\Psi_{X}^{u_1}(h_0) \equiv 1 \mod 2.
\]
Let $\Sigma$ be a closed surface embedded in $X$ which represent the homology class $h_0$. We have the following three cases: (i) $\Sigma \subset X_1$, (ii) $\Sigma \subset X_2$, (iii) $\Sigma \cap Y \not= \emptyset$.

\vspace{2mm}

\noindent
(i) Suppose that $\Sigma \subset X_1$. By Theorem \ref{thm gluing 1}, we have
\[
\Psi_{X}^{u_1}(h_0) = < \Psi_{X_1}^{u_1}([\Sigma]), \Psi_{X_2}^{u_1}>.
\]
It follows from (\ref{eq moduli X_1}) that the relative invariant $\Psi_{X_1}^{u_1}([\Sigma])$ lives in $I_{-1}(Y)$. By Proposition \ref{prop 8N+1}, we have $\Psi_{X_1}^{u_1} ([\Sigma]) = 0$, and hence $\Psi_{X}^{u_1}(h_0) \equiv 0 \mod 2$. This is a contradiction.

\vspace{2mm}

\noindent

\noindent
(ii)
Suppose that $\Sigma \subset X_2$. Then
\[
\Psi_{X}^{u_1}(h_0) = < \Psi_{X_1}^{u_1}, \Psi_{X_2}^{u_1}([\Sigma]) >
\ \text{and} \ 
\Psi_{X_1}^{u_1} \in I_0(Y).
\] 
Since $I_0(Y) = 0$, we have a contradiction.

\vspace{2mm}

\noindent
(iii)
Suppose that $\Sigma \cap Y \not= \emptyset$. We may assume that the intersection $\Sigma \cap Y$ is transverse and the number of  connected components of $\Sigma \cap Y$ is $1$. (If the number of connected components is larger than $1$, join the connected components of $\Sigma \cap Y$ in $Y$ by thin tubes without change of the homology class $[\Sigma] \in H_2(X;\Z)$. )
Thus we can suppose $\gamma := \Sigma \cap Y$ is diffeomorphic to $S^1$. Put $\Sigma_1 := X_1 \cap \Sigma$ and $\Sigma_2 := \Sigma \cap X_2$. Suppose that $\Sigma_1$ ( and hence $\Sigma_2$) satisfies (\ref{eq c1 1}) or (\ref{eq c1 2}). Then by Theorem \ref{thm gluing 2} and (\ref{eq moduli X_1}) we have
\[
\Psi_{X}^{u_1}(h_0) = < \Psi_{X_1}^{u_1}([\Sigma_1]), \Psi_{X_2}^{u_1}([\Sigma_2]) >,
\ \text{and} \ 
\Psi_{X_1}^{u_1}([\Sigma_1]) \in I_0(Y;\gamma).
\]
We show that $\Psi_{X_1}^{u_1}([\Sigma_1]) = 0$ in $I_0(Y;\gamma)$.
Recall that the complex for $I_*(Y;\gamma)$ is as follows:
\begin{equation} \label{eq chain I gamma}
\xymatrix@R=6pt{
C_{1}(Y;\gamma) \ar@{=}[d] \ar[r]^{ \partial_{\gamma} } & C_{0}(Y;\gamma) \ar@{=}[d] \ar[r]^{ \partial_{\gamma} } & C_{-1}(Y;\gamma) \ar@{=}[d] \\
CF_2(Y) \ar[r]^{0} \ar[rdd] \ar[rdddd] & CF_{0}(Y) \ar[r]^{\cong} \ar[rdd] & CF_{-2}(Y) \\
\oplus & \oplus & \oplus \\
CF_{0}(Y) \ar[r]^{\cong} & CF_{-2}(Y) \ar[r]^{0} & CF_{-4}(Y) \\ 
                 & \oplus            &            \\
                 & \Z_2 <[\theta] > \ar[ruu] & 
}
\end{equation}
In the proof of Proposition \ref{prop 8N+1}, we proved that the map $CF_{0}(Y) \rightarrow CF_{-2}(Y)$ is isomorphic and the maps $CF_2(Y) \rightarrow CF_0(Y)$ and $CF_{-2}(Y) \rightarrow CF_{-4}(Y)$ are trivial. (Note that $C_i(Y) = CF_{2i}(Y)$.)

Let $\psi = \psi_{X_1}^{u_1}([\Sigma_1])$ be the element of $C_0(Y;\gamma)$ which represents the class $\Psi_{X_1}^{u_1}([\Sigma_1])$. (See Subsection \ref{ss rel inv} for the definition of $\psi_{X_1}^{u_1}([\Sigma_1])$.) We can write
\[
\psi = \psi_{0} +  \psi_{-2} + n [\theta],
\]
where $\psi_0$ and $\psi_{-2}$ are elements in $CF_0(Y)$ and $CF_{-2}(Y)$ respectively, and $n = <\psi, [\theta]> \in \Z_2$. 
Assume that $\psi_0 \not= 0$, then from the diagram (\ref{eq chain I gamma}), we can see that the $CF_{-2}(Y)$-component of the image $\partial_{\gamma}(\psi)$ is non-trivial. This is a contradiction, since $\psi$ is a cycle. Thus we can write
\[
\psi = \psi_{-2} + n [\theta].
\]
We will prove that $i_{\theta} \equiv 0 \mod 2$ later. (See Lemma \ref{lem i theta} below. Recall that $i_{\theta}$ is the index of Dirac operator twisted by a connection over $\hat{X}_1$ with trivial limit.) By (\ref{eq psi theta}), we have $n \equiv 0 \mod 2$. Hence $\psi$ is included in $CF_{-2}(Y)$. Since $CF_{0}(Y) \rightarrow CF_{-2}(Y)$ is isomorphic, $\psi$ is in the image of the map. We denote the inverse image by $\psi_{0}' \in CF_{0}(Y)$. From the above diagram (\ref{eq chain I gamma}) we can see $\partial_{\gamma}(\psi_0') = \psi$. Therefore $\Psi_{X_1}^{u_1}([\Sigma_1]) = 0$ in $I_0(Y;\gamma)$, and we have a contradiction again.

If $\Sigma_1$ does not satisfy both of (\ref{eq c1 1}) and (\ref{eq c1 2}), then
\begin{eqnarray}
\label{eq c1 3}
&&< c_1(Q_1;\varphi);[\Sigma_1]> \equiv 0 \mod 2 \ 
\text{and} \ 
[\gamma] \not= 0 \ \text{in} \ H_1(Y;\Z), \ \text{or} \\
\label{eq c1 4}
&&< c_1(Q_2;\varphi), [\Sigma_1]> \equiv 1 \mod 2 \ 
\text{and} \ 
[\gamma] = 0 \ \text{in} \ H_1(Y;\Z).
\end{eqnarray}

Assume that (\ref{eq c1 3}) holds. Since  $\Psi_{X_1}^{u_1}$ is a homomorphism, $\Psi_{X_1}^{u_1}(ph_0)$ is also non-trivial. Here $p = 16 N + 1$. Let $\Sigma'$ be an embedded surface in $X$ representing $ph_0$. Then we can easily see that $\Sigma_1' = X_1 \cap \Sigma'$ and $\Sigma_2' = X_2 \cap \Sigma'$ satisfy (\ref{eq c1 2}). Assume that (\ref{eq c1 4}) holds. In this case, we consider the class $h_0 + 2h_1$. Here $h_1 \in H_2(X;\Z)$ is defined as follows. Fix a loop $\gamma_1$ in $Y$ such that the class $[\gamma_1] \in H_1(Y;\Z)$ is a generator. Take relative homology classes $[\Sigma_1''] \in H_2(X_1, Y;\Z)$,  $[\Sigma_2''] \in H_2(X_2,Y;\Z)$ such that $\partial \Sigma_1'' = \gamma_1$, $\partial \Sigma_2'' = \gamma_1$. We can see that there are such surfaces from the exact sequences
\[
H_2(X_i,Y) \longrightarrow H_1(Y) \longrightarrow H_1(X_i) = 0.
\]
Put $h_1 = 2 [\Sigma''_1 \cup \Sigma_2'']$. Then $h_1 \cdot h_1 \equiv 0 \mod 2$ and $\Psi_{X}^{u_1}(h_0 + 2h_1) \equiv 1 \mod 2$. Take a surface $\Sigma'''$ in $X$ which represents the class $h_0 + 2h_1$, and put $\Sigma_1''' := X_1 \cap \Sigma'''$, $\Sigma_2''' = X_2 \cap \Sigma'''$. We can see that $[\Sigma_1''']$ and $[\Sigma_2''']$ satisfy (\ref{eq c1 1}). The same argument as above gives a contradiction in each case.

\vspace{2mm}

It remains to prove the following.

\begin{lem} \label{lem i theta}
Put $Y=L(p,2)$ and let $X_1$ be as in Theorem \ref{thm dec}. Here $p$ is a positive integer of the form $8N + 1$.
Choose a cohomology class $\alpha \in H^2(X_1;\Z)$ as in Lemma \ref{lem alpha} and a $U(2)$-bundle $Q_1$ over $X_1$ with $c_1 = \alpha$. Denote the extension of $Q_1$ to the bundle over $\hat{X}_1 = X_1 \cup (Y \times [0, \infty))$ by $\hat{Q}_1$ as usual. For a connection $\hat{A}_1$ on $\hat{Q}_1$ with trivial limit $\eta$, we have
\[
\ind \Dir_{ \hat{A}_1} \equiv N \mod 2.
\]
Hence if $p$ is of the form $16N + 1$, the index is even.
\end{lem}

To prove this, take a spin 4-manifold $X'$ with boundary $-Y$. For the trivial connection $\eta_{\hat{X}'}$ on the trivial $U(2)$-bundle over $\hat{X}' = X' \cup \big( (-Y) \times [0,\infty) \big)$, we have
\[
\ind \Dir_{\eta_{\hat{X}'}} \equiv 0 \mod 2
\]
since $\Dir_{\eta_{ \hat{X}'}}$ is the direct sum of two copies of the Dirac operator associated with a spin structure. Hence we have
\[
\ind \Dir_{ \hat{A}_1} \equiv 
\ind \Dir_{ \hat{A}_1} + \ind \Dir_{\eta_{\hat{X}'}} \mod 2.
\]
By the additivity of the index, the right hand side is equal to the index $\ind \Dir_{A''}$ of the Dirac operator over $X'' = X_1 \cup_{Y} X'$. Here $A''$ is the connection obtained by gluing $\hat{A}_1$ and $\theta_{\hat{X}'}$. Denote by $Q''$ the $U(2)$-bundle over $X''$ obtained from $Q_1$ and the trivial bundle over $X'$. By Atiyah-Singer's index theorem, 
\[
\begin{split}
\ind \Dir_{A''} 
&= \frac{c_1(Q'')^2 - \sign (X'')}{8} \\
&= \frac{\alpha^2 - \sign(X_1) - \sign(X')}{8} \\
&= \frac{\alpha^2 - 1 - \sign(X')}{8}.
\end{split}
\]
Here $\sign(X'')$, $\sign(X_1)$ and $\sign (X')$ are the signatures of $X''$, $X_1$ and $X'$ respectively, and we have used the additivity of the signatures. We showed that $\alpha^2  = p \ (=8N+1) $ in the proof of Lemma \ref{lem X dec 2}. Hence we have
\[
\ind \Dir_{A''} 
= \frac{p - 1}{8} - \frac{\sign (X')}{8}
= N - \frac{\sign(X')}{8}.
\]
The proof is reduced to showing that $\sign (X') \equiv 0 \mod 16$. We can see this, using Casson-Walker invariant \cite{Wal}. Let $\lambda(M)$ be Casson-Walker invariant of a closed, oriented 3-manifold $M$. Then we have
\[
\sign (X') \equiv 4p^2 \lambda(-L(p,2)) \mod 16.
\]
See \cite[(6.5) Proposition]{Wal}. Moreover we have
\[
\begin{split}
\lambda(-L(p,2)) 
&= 
-\lambda(L(p,2)) 
=
\sum_{k=1}^{p-1} 
\left( \left( \frac{k}{p} \right) \right) 
\left( \left( \frac{2k}{p} \right) \right), \\
((x)) 
&= 
\left\{
\begin{array}{cl}
0 & \text{if $x \in \Z$}, \\
x - [x] - \frac{1}{2} & \text{otherwise}.
\end{array}
\right.
\end{split}
\]
See \cite[(6.3) Proposition]{Wal}. A calculation shows that
\[
\sign(X')
=\frac{p(p-1)(p-5)}{6}
=\frac{16N(2N-1)(8N+1)}{3}.
\]
Since $N(2N-1)(8N+1) \equiv 0 \mod 3$, we have
\[
\sign(X') \equiv 0 \mod 16
\]
as required.

{\it Acknowledgements}. 
The author would like to thank Mikio Furuta, Yuichi Yamada and Kouichi Yasui for useful discussion.


\end{document}